\documentclass[a4paper, 11pt]{article}
\usepackage{amssymb, amsmath,latexsym,amsfonts,amsbsy, amsthm,mathtools,graphicx,CJKutf8,CJKnumb,CJKulem,color,geometry}
\usepackage{verbatim}
\usepackage{hyperref}
\usepackage{newtxtext, newtxmath}
\usepackage[misc]{ifsym}
\usepackage[marginal]{footmisc}
\usepackage[titletoc,title]{appendix}
\usepackage{cite}
\def\l{\left}
\def\r{\right}
\def\f{\frac}

\def\az{\alpha}

\def\az{\alpha}

\def\ez{\epsilon}
\def\bz{\beta}
\def\dz{\delta}
\def\gz{\gamma}

\def\sz{\sigma}

\def\beq{\begin{equation}}
\def\eeq{\end{equation}}
\def\be{\begin{equation*}}
\def\ee{\end{equation*}}
\def\beqn{\begin{eqnarray}}
\def\eeqn{\end{eqnarray}}
\def\ben{\begin{eqnarray*}}
\def\een{\end{eqnarray*}}

\theoremstyle{plain}
\theoremstyle{plain}\newtheorem{thm}{Theorem}[section]
\theoremstyle{plain}\newtheorem{prop}{Proposition}[section]
\theoremstyle{plain}
\theoremstyle{plain}\newtheorem{lem}{Lemma}[section]
\theoremstyle{plain}

\numberwithin{equation}{section}

\makeatletter
\def\thanks#1{\protected@xdef\@thanks{\@thanks
        \protect\footnotetext{#1}}}
\makeatother

\begin{document}

\title{Global solutions to quasilinear wave-Klein-Gordon systems in two space dimensions}
\author{Qian Zhang}

\date{}

\maketitle

\noindent {\bf{Abstract}}\ \ In this paper we prove global existence and global behavior of solutions to quasilinear wave-Klein-Gordon systems in $\mathbb{R}^{1+2}$ with quadratic nonlinearities satisfying the null condition. We consider small, regular and compactly supported initial data, and prove global existence, pointwise decay estimates and linear scattering  for the solutions.
\bigskip


\bigskip


\section{Introduction}\label{sI}

We consider the following quasilinear wave-Klein-Gordon system in $\mathbb{R}^{1+2}$
\beq\label{qEquv}
\l\{\begin{array}{rcl}
-\Box u&=&F_u:=P_1^{\gz\az\bz}\partial_\gz v\partial_\az\partial_\bz v+P_2^{\gz\az\bz}\partial_\gz u\partial_\az\partial_\bz v=N_1(v,v)+N_2(u,v),\\
-\Box v+v&=&F_v:=P_1^{\gz\az\bz}\partial_\gz v\partial_\az\partial_\bz u+P_2^{\gz\az\bz}\partial_\gz u\partial_\az\partial_\bz u=N_1(v,u)+N_2(u,u)
\end{array}\r.
\eeq
with compactly supported initial data prescribed on the time slice $t=t_0=2$:
\beq\label{qini}
(u,\partial_tu,v,\partial_tv)|_{t=t_0}=(u_0,u_1,v_0,v_1).
\eeq
Here $\gz,\az,\bz\in\{0,1,2\}$ and we denote $N_i(w,z)=P_i^{\gz\az\bz}\partial_\gz w\partial_\az\partial_\bz z$ for any smooth functions $w,z$, for $i=1,2$. Einstein summation convention over repeated upper and lower indices is adopted throughout the paper. In the above, $P_1^{\gz\az\bz},P_2^{\gz\az\bz}$ are constants satisfying the standard null condition, that is, 
\beq\label{nullcond}
P_i^{\gz\az\bz}\xi_\gz\xi_\az\xi_\bz=0\quad\mathrm{for\; all\;\;} \xi_0^2=\xi_1^2+\xi_2^2,\quad i=1,2,
\eeq
and in addition, the symmetry condition 
\be
P_i^{\gz\az\bz}=P_i^{\gz\bz\az},\quad i=1,2.
\ee
As usual, $\Box=g^{\az\bz}\partial_\az\partial_\bz=-\partial_0^2+\partial_1^2+\partial_2^2$ denotes the wave operator, where $\partial_0=\partial_t,\partial_a=\partial_{x_a}$ for $a=1,2$, $g=(g_{\az\bz})=\mathrm{diag}(-1,1,1)$ denotes the Minkowski metric in $\mathbb{R}^{1+2}$, $(g^{\az\bz})$ denotes the inverse matrix of $(g_{\az\bz})$. Without loss of generality, the initial data $(u_0,u_1,v_0,v_1)$ are assumed to be supported in the unit ball $\{x: |x|<1\}$, hence the solution is supported within the region $\{(t,x): t\ge 2,|x|<t-1\}$.

Throughout this paper, Greek letters $\gz,\az,\bz,\dots\in\{0,1,2\}$ represent spacetime indices and Latin letters $a,b,c,\dots\in\{1,2\}$ are used for space indices. For any two quantities $A,B\ge 0$, we write $A\lesssim B$ if $A\le CB$ for some unimportant constant $C>0$. We write $A\sim B$ if $A\lesssim B$ and $B\lesssim A$. We write $A\ll B$ if $A\le CB$ for some constant $C>0$ sufficiently small.

Let us first review some works on global well-posedness for nonlinear wave equations. For general nonlinear wave equations in $\mathbb{R}^{1+3}$ with quadratic nonlinearities, the local solution may blow up in finite time; see John \cite{J81} for example. Klainerman \cite{K86} and Christodoulou \cite{Ch} independently proved that nonlinear wave equations in $\mathbb{R}^{1+3}$ with small data and nonlinearities satisfying the null condition admit global-in-time solutions. Using the "ghost weight" energy estimates, Alinhac \cite{Al1} established global existence for quasilinear wave equations in $\mathbb{R}^{1+2}$ with quadratic null nonlinearities, for small, smooth and compactly supported initial data. In \cite{HY}, Hou-Yin removed the compactness assumption on the support of initial data in \cite{Al1}. A similar result for $2D$ fully nonlinear wave equations under the null condition was obtained by Cai, Lei, and Masmoudi \cite{CLM}. In \cite{Al1}, the top order energy of the solution grows polynomially in time. A question, known as the "blowup-at-infinity" conjecture \cite{Al02,Al03,Al10}, is whether this growth is a true phenomenon. This was solved by Dong-LeFloch-Lei \cite{DLL} and Li \cite{Ld} independently, where it was shown that the top order energy of quasilinear wave equations in $\mathbb{R}^{1+2}$ with quadratic null nonlinearities is uniformly bounded in time.  

We now mention some related works on coupled wave-Klein-Gordon systems. In \cite{Geo}, Georgiev proved global existence of small solutions for coupled systems of nonlinear wave and Klein–Gordon equations with strong null condition in $\mathbb{R}^{1+3}$, which was improved by Katayama \cite{Kat} to more general nonlinearities. LeFloch-Ma \cite{LM16}, Wang \cite{Wa} and Ionescu-Pausader \cite{Ion} studied the wave-Klein-Gordon system in $\mathbb{R}^{1+3}$ as a model for the full Einstein–Klein–Gordon system. In lower space dimensions, Ma \cite{M17,M17'} studied global existence for a quasilinear diagonalized wave-Klein–Gordon system in $\mathbb{R}^{1+2}$, and then extended the result \cite{M17'} to more types of nonlinearities in \cite{M19,M21}. In \cite{St18}, Stingo proved global existence for a quasilinear wave-Klein-Gordon system in $\mathbb{R}^{1+2}$ with $Q_0$ type nonlinearities (here $Q_0(w,z)=\partial_\az w\partial^\az z$ for any functions $w,z$), when initial data are small, smooth and mildly decay at infinity. Ifrim-Stingo \cite{IS} established almost global existence for quasilinear wave-Klein-Gordon systems in $\mathbb{R}^{1+2}$ with quadratic null nonlinearities. There also exist many other results on nonlinear wave and Klein-Gordon equations as well as their coupled systems; see for instance \cite{Al2,Lind,DFX,Zha,OTT,FWY,D21,KWY,Do21,DWd,DuM,DW20}.

Inspired by the works \cite{IS,St18}, the goal of this paper is to prove global existence and asymptotic behavior for the solution to \eqref{qEquv}-\eqref{qini} under the null condition (i.e., \eqref{nullcond}), for small, regular and compactly supported initial data. We use the hyperboloidal method which is due to Klainerman \cite{K2,K93} and H\"ormander \cite{Ho}.

Before we present our main results, some notations are made as follows. 

Let $l\in\mathbb{N}$. We denote ${\bf H}^l(\mathbb{R}^2):=H^{l+1}(\mathbb{R}^2)\times H^l(\mathbb{R}^2)$ and ${\bf{H}}_{l}(\mathbb{R}^2):=\big(\dot{H}^{l+1}(\mathbb{R}^2)\cap\dot{H}^1(\mathbb{R}^2)\big)\times H^l(\mathbb{R}^2)$, where $H^k(\mathbb{R}^2), \dot{H}^{k}(\mathbb{R}^2), k\in\mathbb{N}$ denote the Sobolev spaces and homogeneous Sobolev spaces respectively. We denote 
\beq\label{defX_l}
\mathcal{X}_{l}(\mathbb{R}^2):={\bf{H}}_{l}(\mathbb{R}^2)\times {\bf{H}}^{l}(\mathbb{R}^2)=\big(\dot{H}^{l+1}(\mathbb{R}^2)\cap\dot{H}^1(\mathbb{R}^2)\big)\times H^l(\mathbb{R}^2)\times H^{l+1}(\mathbb{R}^2)\times H^l(\mathbb{R}^2).
\eeq

We are now ready to state the main results of this paper.
 
\begin{thm}\label{thm1}
Let $N\ge 14$ be an integer. Consider the quasilinear wave-Klein-Gordon system \eqref{qEquv} with initial data $(u_0,u_1,v_0,v_1)$ on the time slice $t=t_0=2$ supported in the ball $\{x: |x|<1\}$. Then for any $\dz>0$, there exists $\ez_0>0$ such that, for all $0<\ez\le\ez_0$ and all initial data satisfying
\be
\|u_0\|_{H^{N+1}(\mathbb{R}^2)}+\|u_1\|_{H^{N}(\mathbb{R}^2)}+\|v_0\|_{H^{N+1}(\mathbb{R}^2)}+\|v_1\|_{H^{N}(\mathbb{R}^2)}\le\ez,
\ee
the Cauchy problem \eqref{qEquv}-\eqref{qini} admits a global-in-time solution $(u,v)$, which satisfies the following pointwise decay estimates
\be
|v(t,x)|\lesssim t^{-1},\quad\quad |u(t,x)|\lesssim t^{-1/2+\dz},\quad\quad |\partial u(t,x)|\lesssim t^{-1/2}.
\ee
Moreover, the solution $(u,v)$ scatters to a free solution in $\mathcal{X}_{N-5}(\mathbb{R}^2)$ (see  \eqref{defX_l}), i.e., there exists $(u^*_0,u^*_1,v^*_0,v^*_1)\in\mathcal{X}_{N-5}(\mathbb{R}^2)$ such that
\ben
\lim_{t\to+\infty}\|(u,\partial_tu,v,\partial_tv)-(u^*,\partial_tu^*,v^*,\partial_tv^*)\|_{\mathcal{X}_{N-5}(\mathbb{R}^2)}=0,
\een
where $(u^*,v^*)$ is the solution to the $2D$ linear homogeneous wave-Klein-Gordon system with the initial data $(u^*_0,u^*_1,v^*_0,v^*_1)$.
\end{thm}

\noindent${\bf{Difficulties\;and\;key\;ideas.}}$ We follow \cite{K2,Ho} and use hyperboloids $\mathcal{H}_s=\{(t,x): t^2=s^2+|x|^2\}$ ($s\ge s_0=2$) to foliate the spacetime. Energy estimates are derived along these hyperboloids and integration is with respect to the hyperbolic time $s=\sqrt{t^2-|x|^2}$ (instead of $t$). The main advantage of this approach is that we can make use of the $(t-|x|)$ decay. To prove the global existence result in Theorem \ref{thm1}, the main challenges include the following: $i)$ The nondiagonalizable (i.e., $F_u$ contains $\partial\partial v$ and $F_v$ contains $\partial\partial u$) structure of the nonlinearities brings difficulty in deriving an inequality for the top order energy of the solution $(u,v)$ to \eqref{qEquv}-\eqref{qini}; $ii)$ The slow decay nature of quadratic nonlinearities causes trouble in closing the energy estimate. To close the estimate of top order energy (which is expected to have a small growth), we need to obtain sharp decay estimates for the solution. For this, we need a uniform (in time) bound for lower order energy of the solution. To overcome these difficulties, we adopt some novel ideas as stated below.

First, to derive an inequality for the energy of the solution $(u,v)$ up to the top order, we combine both equations in \eqref{qEquv} rather than dealing with single equation. Precisely, by acting the vector field $\Gamma^I$ ($\Gamma\in\{\partial_\az,L_a\}$ and $I$ is a multi-index, where $\partial_\az, L_a$ denote the translations and Lorentz boosts respectively) on both sides of each equation in \eqref{qEquv}, and applying Lemma \ref{lemHo}, we obtain
\beq\label{Eqgzu'}
-\Box \Gamma^Iu=\Gamma^IF_u=\sum_{I_1+I_2\le I}\l\{N_{1,I;I_1,I_2}(\Gamma^{I_1}v,\Gamma^{I_2}v)+N_{2,I;I_1,I_2}(\Gamma^{I_1}u,\Gamma^{I_2}v)\r\},
\eeq
\beq\label{Eqgzv'}
-\Box \Gamma^Iv+\Gamma^Iv=\Gamma^IF_v=\sum_{I_1+I_2\le I}\l\{N_{1,I;I_1,I_2}(\Gamma^{I_1}v,\Gamma^{I_2}u)+N_{2,I;I_1,I_2}(\Gamma^{I_1}u,\Gamma^{I_2}u)\r\},
\eeq
where for $i=1,2$ and any sufficiently smooth functions $w,z$, $N_{i,I;I_1,I_2}(w,z)=P_{\!i,I;I_1,I_2}^{\gz\az\bz}\partial_\gz w\partial_\az\partial_\bz z$ with $P_{\!i,I;I_1,I_2}^{\gz\az\bz}$ satisfying the null condition, and $N_{i,I;I_1,I_2}(w,z)=N_i(w,z)$ when $I_1+I_2=I$. Multiplying \eqref{Eqgzu'} and \eqref{Eqgzv'} by $\partial_t\Gamma^Iu$ and $\partial_t\Gamma^Iv$ respectively, we have
\beq\label{Equdtu'}
\f{1}{2}\partial_t(|\partial\Gamma^Iu|^2)-\partial_a(\partial^a\Gamma^Iu\partial_t\Gamma^Iu)=F_u^{I,low}\partial_t\Gamma^Iu+\big(N_1(v,\Gamma^Iv)+N_2(u,\Gamma^Iv)\big)\partial_t\Gamma^Iu,
\eeq
\beq\label{Eqvdtv'}
\f{1}{2}\partial_t(|\partial\Gamma^Iv|^2+|\Gamma^Iv|^2)-\partial_a(\partial^a\Gamma^Iv\partial_t\Gamma^Iv)=F_v^{I,low}\partial_t\Gamma^Iv+\big(N_1(v,\Gamma^Iu)+N_2(u,\Gamma^Iu)\big)\partial_t\Gamma^Iv
\eeq
(see \eqref{dubdu} for the definitions of $|\partial\Gamma^Iu|$ and $|\partial\Gamma^Iv|$), where 
\ben
F_u^{I,low}:&=&\sum_{I_1+I_2\le I,I_2|\le |I|-1}\big(N_{1,I;I_1,I_2}(\Gamma^{I_1}v,\Gamma^{I_2}v)+N_{2,I;I_1,I_2}(\Gamma^{I_1}u,\Gamma^{I_2}v)\big),\\
F_v^{I,low}:&=&\sum_{I_1+I_2\le I,|I_2|\le |I|-1}\big(N_{1,I;I_1,I_2}(\Gamma^{I_1}v,\Gamma^{I_2}u)+N_{2,I;I_1,I_2}(\Gamma^{I_1}u,\Gamma^{I_2}u)\big). 
\een
By adding \eqref{Equdtu'} and \eqref{Eqvdtv'}, and a careful calculation in the terms "$N_1(v,\Gamma^Iv)\partial_t\Gamma^Iu+N_1(v,\Gamma^Iu)\partial_t\Gamma^Iv$" and "$N_2(u,\Gamma^Iv)\partial_t\Gamma^Iu+N_2(u,\Gamma^Iu)\partial_t\Gamma^Iv$", we obtain an equality where the terms $\partial\partial\Gamma^Iu,\partial\partial\Gamma^Iv$ vanish on the right hand side.     After integrating over the region limited by two hyperboloids, we obtain an energy equality for $(\Gamma^Iu,\Gamma^Iv)$.
 
Next, to close the energy estimate up to the top order, we need to gain a uniform in time bound for lower order energy of the solution $(u,v)$. This is achieved by performing nonlinear transformations for both $u$ and $v$. To bound the lower order energy of $v$, we let $\tilde{v}=v-N_2(u,u)$. Then $\tilde{v}$ solves  
\beq\label{boxtv}
-\Box\tilde{v}+\tilde{v}=N_1(v,u)+N_2(\partial u, \partial u)+"\mathrm{good\;terms}"
\eeq
(here we omit the constant coefficients of each term), where we refer to "good terms" as "cubic terms involving $v$ or derivatives of $v$". Now in the term $N_2(\partial u,\partial u)$ there is one "$\partial$" hitting each $u$. Hence, by the null condition and extra decay of Hessian of the wave component $u$ (see Lemmas \ref{lemnf} and \ref{lemddu}), the terms $N_1(v,u)$ and $N_2(\partial u, \partial u)$ appearing on the right hand side of \eqref{boxtv} have sufficient decay rates.

To bound the lower order energy of $u$, the main difficulty is to take care of the term $N_1(v,v)=P_1^{\gz\az\bz}\partial_\gz v\partial_\az\partial_\bz v$ in \eqref{qEquv}. For this, we carefully exploit the structure of the nonlinearities, and discover new nonlinear transformations leading to faster decay nonlinearities. Precisely, we substitute  the $"v"$ appearing in $"\partial_\gz v"$ in $N_1(v,v)$ by $\Box v+F_v$ (here we use the second equation in \eqref{qEquv}), and arrive at
\beq\label{boxu}
-\Box u=N_1(-\partial_t\partial_t v+\partial_a\partial^av,v)+N_2(u,v)+"\mathrm{good\;terms}".
\eeq
To deal with the term $N_1(-\partial_t\partial_t v+\partial_a\partial^av,v)$ on the right hand side of \eqref{boxu}, we compute $(-\Box)N_1(v,v)$ and make a subtle cancellation. On the other hand, the term $N_2(u,v)$ on the right hand side of \eqref{boxu} is cancelled by conducting the transformation $u+N_2(u,v)$. Combining both terms, our final transformation is $\tilde{u}:=u+\f{1}{4}N_1(v,v)+N_2(u,v)$. The nonlinearity $F_{\tilde{u}}:=-\Box\tilde{u}$ for this new function has sufficient decay rate. Hence we can close the bootstrap and obtain the global existence result.

Finally, to show the scattering result in Theorem \ref{thm1}, we need to estimate the $L^2_x$ norms of the nonlinearities on flat time slices. However, the hyperboloidal method we use only provides estimates on hyperboloids. Hence, we prove a technical lemma which gives a sufficient condition on the $L^2$-type norms of the nonlinearities on hyperboloids, for the linear scattering of the solution.

The organization of this paper is as follows. In Section \ref{sP}, we introduce some notations, and state energy and Sobolev inequalities on hyperboloids, and estimates of null forms. In Section \ref{sE}, we provide the main ingredients in proving Theorem \ref{thm1}, including an energy equality and some nonlinear transformations. The last two sections are devoted to the complete proof of Theorem \ref{thm1}. Precisely, we prove the global existence and the linear scattering results in Theorem \ref{thm1} in Sections \ref{sM} and \ref{sSca} respectively.

\section{Preliminaries}\label{sP}

\subsection{Notations}\label{sno}
We work in the $(1+2)$ dimensional spacetime $\mathbb{R}^{1+2}$ with Minkowski metric $g=(-1,1,1)$, which is used to raise or lower indices. We denote a point in $\mathbb{R}^{1+2}$ by $(t,x)=(x_0,x_1,x_2)$ with $t=x_0,x=(x_1,x_2),x^a=x_a,a=1,2$, and its spacial radius is denoted by $r:=|x|=\sqrt{x_1^2+x_2^2}$. Following Klainerman's vector field method \cite{K86}, we introduce the following vector fields:
\begin{itemize}
\item[(i)] Translations: $\partial_\az:=\partial_{x_\az}$,  $\az\in\{0,1,2\}$.
\item[(ii)] Lorentz boosts: $L_a:=x_a\partial_t+t\partial_a$,  $a\in\{1,2\}$.
\item[(iii)] Rotation: $\Omega_{12}:=x_1\partial_2-x_2\partial_1$.
\item[(iv)] Scaling: $L_0=t\partial_t+x^a\partial_a$.
\end{itemize}
For any operators $A$ and $B$, the commutator $[A,B]$ is defined as
$$[A,B]:=AB-BA.$$
For simplicity, we denote $\sum_\az=\sum_{\az\in\{0,1,2\}}$ and similarly for $\sum_\bz,\sum_\gz$, while $\sum_a=\sum_{a\in\{1,2\}}$ and similar for $\sum_b,\sum_c$.

We restrict our study to functions supported within the spacetime region  
$$\mathcal{K}:=\{(t,x): t\ge 2, r<t-1\},$$ 
which is the light cone with vertex $(1,0,0)$. We use $s\ge s_0=2$ to denote the hyperbolic time, and hyperboloids are denoted by
\be
\mathcal{H}_s:=\{(t,x): t^2-r^2=s^2\}.
\ee
We note that for any $(t,x)\in\mathcal{K}\cap\mathcal{H}_s$ with $s\ge s_0=2$, it holds that
\be\label{sC: tles^2}
r\le t-1,\quad\quad s\le t\le s^2.
\ee
For any $s_1\ge s_0=2$, we denote by $\mathcal{K}_{[s_0,s_1]}:=\bigcup_{s_0\le s\le s_1}\mathcal{K}\cap\mathcal{H}_s$ the subsets of $\mathcal{K}$ limited by the hyperboloids $\mathcal{H}_{s_0}$ and $\mathcal{H}_{s_1}$. We follow \cite{LM14} and introduce the hyperboloidal frame, which is defined by
\beq\label{defbdaz}
\bar{\partial}_0=\partial_s=\f{s}{t}\partial_t,\quad\quad\bar{\partial}_a=\f{L_a}{t}=\f{x_a}{t}\partial_t+\partial_a,\quad a=1,2.
\eeq
We also make use of the semi-hyperboloidal frame 
\be
\underline{\partial}_0:=\partial_t,\quad\quad\underline{\partial}_a=\bar{\partial}_a,\quad a=1,2.
\ee
For any sufficiently smooth function $u$, we denote for simplicity
\beq\label{dubdu}
|\partial u|=\Bigg(\sum_{\az}|\partial_\az u|^2\Bigg)^{1/2},\quad\quad\quad\quad|\bar{\partial}u|=\Bigg(\sum_{\az}|\bar{\partial}_\az u|^2\Bigg)^{1/2}.
\eeq
Given a sufficiently nice function $u$ supported in $\mathcal{K}$, its $L^p$ norms on the hyperboloids $\mathcal{H}_s$ ($s\ge s_0=2$) are defined by
\be\label{L^2fHtau}
\|u\|_{L^p_f(\mathcal{H}_s)}^p=\int_{\mathcal{H}_s}|u(t,x)|^p{\rm{d}}x:=\int_{\mathbb{R}^2}|u(\sqrt{s^2+|x|^2},x)|^p{\rm{d}}x,\quad 1\le p<\infty.
\ee

We denote the ordered set
\be
\{\Gamma_k\}_{k=0}^{4}:=\{\partial_0,\partial_1,\partial_2, L_1,L_2\}.
\ee
For any multi-index $I=(i_0,i_1,i_2,i_3,i_4)\in\mathbb{N}^{5}$ of length $|I|=\sum_{k=0}^{4}i_k$, we denote
\be
\Gamma^I=\prod_{k=0}^{4}\Gamma_k^{i_k},\quad\mathrm{where}\quad\Gamma=(\Gamma_0,\Gamma_1,\Gamma_2,\Gamma_3,\Gamma_4).
\ee
For any multi-indices $I=(i_0,i_1,i_2)\in\mathbb{N}^{3}, J=(j_1,j_2)\in\mathbb{N}^2$, let 
\beq\label{defdI}
\partial^I:=\prod_{k=0}^{2}\partial_k^{i_{k}},\quad\quad\quad\quad L^J:=\prod_{k=1}^{2}L_k^{j_{k}}.
\eeq
For any multi-indices $I=(i_0,i_1,i_2,i_3,i_4), J=(j_0,j_1,j_2,j_3,j_4)\in\mathbb{N}^{5}$, by writing "$I\le J$" we mean that 
\be
i_k\le j_k\quad\quad\mathrm{for\;all\;\;} k=0,1,2,3,4.
\ee

\subsection{Energy and Sobolev inequalities}

The following estimates for commutators will be frequently used in the sequel.

\begin{lem}(See \cite{LM16})\label{sC: DL2.3}
For any sufficiently smooth function $u$ supported in $\mathcal{K}$, and any multi-indices $I, J\in\mathbb{N}^{5},K\in\mathbb{N}^2$, we have
\ben
&&|\Gamma^I\Gamma^Ju-\Gamma^{I+J}u|\lesssim\sum_{|K'|<|I|+|J|}|\Gamma^{K'}u|,\\
&&|L^KL_a\Gamma^Iu|\lesssim\sum_{b}\sum_{|K'|\le|I|+|K|}|L_b\Gamma^{K'}u|,\quad a\in\{1,2\},\\
&&|L^K\partial_\az\Gamma^Iu|\lesssim\sum_{\bz}\sum_{|K'|\le|I|+|K|}|\partial_\bz\Gamma^{K'}u|,\quad \az\in\{0,1,2\},\\
&&|L^K\underline{\partial}_a\Gamma^Iu|\lesssim\sum_{b}\sum_{|K'|\le|I|+|K|}|\underline{\partial}_b\Gamma^{K'}u|,\quad a\in\{1,2\},\\
\een
\begin{proof}
The proof can be found in \cite{LM14,LM16}, and we sketch it below. Noting that
\be
[\partial_\az,L_b]=\dz_{\az 0}\partial_b+\dz_{\az b}\partial_t,\quad\quad [L_a,L_b]=\f{x_a}{t}L_b-\f{x_b}{t}L_a,
\ee
and that $|L_c(x_a/t)|\lesssim 1$, $|\partial_\az(x_a/t)|\lesssim 1$, $|L_c(1/t)|\lesssim1/t$ in $\mathcal{K}$, for $\az\in\{0,1,2\}, a,c\in\{1,2\}$, we obtain the conclusions.
\end{proof}
\end{lem}

Let $m\ge 0$. Following \cite{LM14,M19,Wo}, for a function $u$ supported in $\mathcal{K}$, we define its (natural) energy and conformal energy on the hyperboloids $\mathcal{H}_s$ ($s\ge s_0=2$) by
\beq\label{DefEner}
\begin{split}
\mathcal{E}_m(u,s):=&\ \int_{\mathcal{H}_s}\Big(|\partial_t u|^2+\sum_{a}|\partial_au|^2+2(x^a/t)\partial_tu\partial_au+m^2u^2\Big){\rm{d}}x\\
=&\ \int_{\mathcal{H}_s}\Big(|(s/t)\partial_tu|^2+\sum_{a}|\underline{\partial}_au|^2+m^2u^2\Big){\rm{d}}x,
\end{split}
\eeq
\be\label{DefCon}
\!\!\!\!\!\!\!\!\!\!\!\!\!\!\!\!\!\!\!\!\!\!\!\!\!\!\!\!\!\!\!\!\!\!\!\!\!\!\!\!\!\!\!\!\!\!\!\!\mathcal{E}_{con}(u,s):=\int_{\mathcal{H}_s}\Big((K_0u+u)^2+\sum_{a}(s\bar{\partial}_au)^2\Big){\rm{d}}x
\ee
respectively, where 
\beq\label{K_0}
K_0:=s\partial_s+2x^a\bar{\partial}_a.
\eeq
Using that $\partial_a=\underline{\partial}_a-(x_a/t)\partial_t$ and recalling \eqref{dubdu}, we have
\beq\label{s/tduE}
\|(s/t)|\partial u|\ \!\|_{L^2_f(\mathcal{H}_s)}+\|\ \!|\bar{\partial}u|\ \!\|_{L^2_f(\mathcal{H}_s)}\lesssim[\mathcal{E}_0(u,s)]^{1/2}.
\eeq

\begin{prop}(See \cite{LM14}.)\label{sC: DL5.1}
Let $m\ge 0$. For any sufficiently smooth function $u$ supported in $\mathcal{K}$, and any $s\in[s_0,+\infty)$, we have
\be
[\mathcal{E}_m(u,s)]^{1/2}\lesssim [\mathcal{E}_m(u,s_0)]^{1/2}+\int_{s_0}^s\|-\Box u+m^2u\|_{L^2_f(\mathcal{H}_\tau)}{\rm{d}}\tau.
\ee
\begin{proof}
For completeness, we sketch the proof below. Let $F:=-\Box u+m^2u$. Multiplying on both sides of this equality by $\partial_t u$, we obtain 
\beq\label{s2: d_tubo}
\partial_t\big(|\partial u|^2+m^2u^2\big)-2\partial_a\big(\partial_tu\partial^au\big)=2\partial_tu F,
\eeq
where we recall \eqref{dubdu}. Note that on $\mathcal{H}_{s}$ we have ${\bf n}{\rm{d}}\sz=(1,-x/t){\rm{d}}x$, where ${\bf n}$ and ${\rm{d}}\sz$ denote the upward unit normal and the volume element of $\mathcal{H}_s$ respectively. Integrating \eqref{s2: d_tubo} over the region $\mathcal{K}_{[s_0,s]}$, and using the transformation $(t,x)\to (\tau,x)$, where $\tau=\sqrt{t^2-|x|^2}$, we obtain 
\be
\mathcal{E}_{m}(u,s)=\mathcal{E}_{m}(u,s_0)+2\int_{s_0}^s\int_{\mathcal{H}_\tau}\f{\tau}{t}\partial_tuF{\rm{d}}x{\rm{d}}\tau.
\ee
Differentiating the last equality with respect to $s$ and using H\"older inequality, the conclusion follows.
\end{proof}
\end{prop}

We next recall the conformal energy estimate on hyperboloids, which were proved by Wong \cite{Wo} and also Ma \cite{M19}.

\begin{prop}(See \cite{Wo,M19,DLL}.)\label{sC: DL5.1con}
For any sufficiently smooth function $u$ supported in $\mathcal{K}$, and any $s\in[s_0,+\infty)$, we have
\be
[\mathcal{E}_{con}(u,s)]^{1/2}\lesssim [\mathcal{E}_{con}(u,s_0)]^{1/2}+\int_{s_0}^s\tau\|\Box u\|_{L^2_f(\mathcal{H}_\tau)}{\rm{d}}\tau.
\ee
\begin{proof}
We only sketch the proof. By straightforward computation, we can express the wave operator $-\Box$ in terms of the hyperboloidal frame as follows
\ben
-\Box u=s^{-1}\partial_s\l(s\partial_su+2x^a\bar{\partial}_au+u\r)-\bar{\partial}_a\bar{\partial}^au.
\een
Using the definition of $K_0$ in \eqref{K_0} and by direct calculation, we arrive at
\ben
s\l(K_0u+u\r)(-\Box u)&=&\f{1}{2}\partial_s\l(K_0u+u\r)^2+\f{1}{2}\partial_s\l(s^2\bar{\partial}_au \bar{\partial}^au\r)-s^2\bar{\partial}_a\l(\partial_su\bar{\partial}^au\r)\\
&-&2s\bar{\partial}_b\l(x^a\bar{\partial}_au\bar{\partial}^bu\r)+s\bar{\partial}_a\l(x^a\bar{\partial}_bu\bar{\partial}^bu\r)-s\bar{\partial}_a\l(u\bar{\partial}^au\r)
\een
Integrating the above identity over $\mathcal{K}_{[s_0,s]}$, we obtain the conclusion.
\end{proof}
\end{prop}

The result below is concerned with the weighted $L^2$-type estimate for a function $u$ on hyperboloids.

\begin{prop}(See \cite{M19}.)\label{sC: DL5.2}
For any sufficiently smooth function $u$ supported in $\mathcal{K}$, and any $s\in[s_0,+\infty)$, we have
\be
\l\|\f{s}{t}u\r\|_{L^2_f(\mathcal{H}_s)}\lesssim\l\|\f{s_0}{t}u\r\|_{L^2_f(\mathcal{H}_{s_0})}+\int_{s_0}^s\f{[\mathcal{E}_{con}(u,\tau)]^{1/2}}{\tau}{\rm{d}}\tau.
\ee
\begin{proof}
We only sketch the proof. By direct computation, we have the following differential identity:
\ben
\partial_s\l(\f{s^2}{t^2}u^2\r)+\bar{\partial}_a\l(\f{x^a}{s}\f{s^2}{t^2}u^2\r)=\f{2}{s}\l[\l(K_0u+u\r)\cdot\f{s^2}{t^2}u-\f{x^a}{t}s\bar{\partial}_au\cdot\f{s}{t}u\r].
\een
Integrating the above identity over $\mathcal{K}_{[s_0,s]}$ yields the conclusion.
\end{proof}
\end{prop}

We next present the Sobolev-type inequalities on hyperboloids which have been proved by Klainerman \cite{K2}, H\"ormander \cite{Ho} and LeFloch-Ma \cite{LM14}. We give the version of LeFloch-Ma, in which only the vector fields of Lorentz boosts are used. 

\begin{lem}(See \cite{LM14,LM16}.)\label{sC: DL2.4}
For any sufficiently smooth function $u$ supported in $\mathcal{K}$, and any $s\in[s_0,+\infty)$, we have
\ben
&&\sup_{\mathcal{H}_s}|tu(t,x)|\lesssim \sum_{|J|\le 2}\|L^Ju\|_{L^2_f(\mathcal{H}_s)},\\
&&\sup_{\mathcal{H}_s}|su(t,x)|\lesssim\sum_{|J|\le 2}\l\|\f{s}{t}L^Ju\r\|_{L^2_f(\mathcal{H}_s)}.
\een
\begin{proof}
The proof of the first inequality can be found in \cite{LM16}. For the second one, let $\tilde{u}(t,x)=\f{s}{t}u(t,x)$. By direct calculation, 
\ben
&&\l[L_a,\f{s}{t}\r]=t\partial_a\l(\f{s}{t}\r)+x_a\partial_t\l(\f{s}{t}\r)=-\f{x_as}{t^2},\\
&&\l[L_bL_a,\f{s}{t}\r]=L_b\l[L_a,\f{s}{t}\r]+\l[L_b,\f{s}{t}\r]L_a=\l(\f{2x_ax_b}{t^2}-\dz_{ab}\r)\f{s}{t}-\f{x_bs}{t^2}L_a.
\een
Applying the first inequality to $\tilde{u}$, we obtain the second inequality.
\end{proof}
\end{lem}

Below we state the extra decay for Hessian of the wave component.

\begin{lem}(See \cite{LM14,M21})\label{lemddu}
Let $w$ solve the wave equation
\be
-\Box w=F_w,
\ee
then we have
\be
|\partial\partial w|\lesssim\f{t}{s^2}\sum_{|I|\le 1}|\partial\Gamma^Iw|+\f{t^2}{s^2}|F_w|\quad\mathrm{in\;}\mathcal{K}.
\ee
\begin{proof}
For completeness, we revisit the proof in \cite{LM14}. We write the d'Alembert operator $-\Box$ as
\be
-\Box=\f{(t-r)(t+r)}{t^2}\partial_t\partial_t+\f{x^a}{t^2}\partial_tL_a-\f{1}{t}\partial^aL_a+\f{2}{t}\partial_t-\f{x^a}{t^2}\partial_a,
\ee
which implies 
\beq\label{dtdtw}
|\partial_t\partial_tw|\lesssim\f{1}{t-r}\sum_{|I|\le 1}|\partial\Gamma^Iw|+\f{t}{t-r}|-\Box w|.
\eeq
We also have
\ben
\partial_t\partial_aw&=&\partial_t\big(t^{-1}L_aw-(x_a/t)\partial_tw\big)=t^{-1}\partial_tL_aw-(x_a/t)\partial_t\partial_tw-t^{-1}\partial_aw,\\
\partial_b\partial_aw&=&\partial_b\big(t^{-1}L_aw-(x_a/t)\partial_tw\big)=t^{-1}\partial_bL_aw-(x_a/t)\partial_b\partial_tw-\dz_{ab}t^{-1}\partial_tw.
\een
Hence \eqref{dtdtw} also holds if we replace $\partial_t\partial_tw$ by $\partial_t\partial_aw$ and $\partial_b\partial_aw$. Noting that $(t-r)\sim s^2/t$ in $\mathcal{K}$, we obtain the conclusion. 
\end{proof}
\end{lem}

\subsection{Estimates of null forms}

The lemma below provides estimates of the null forms appearing in \eqref{qEquv}, whose proof can be found in \cite{LM14,DLL}. For completeness, we revisit the proof and present more conclusions.

\begin{lem}\label{lemnf}
Let $P^{\gz\az\bz}$ satisfy the null condition. Then the following statements hold:
\begin{itemize}
\item[i)] For any sufficiently smooth functions $u,v$ supported in $\mathcal{K}$, we have
\be
P^{\gz\az\bz}\partial_\gz v\partial_\az\partial_\bz u=g(t,x)\partial_tv\partial_t\partial_tu+B(v,u),
\ee
where $g(t,x)$ satisfies
\be
|g(t,x)|\lesssim\f{s^2}{t^2},\quad\quad|\Gamma^I(g(t,x))|\lesssim 1\quad\mathrm{for\;any\;multi\!-\!index\;}I,
\ee
and $B(v,u)$ can be written as
\be
B(v,u)=\f{1}{t}\sum_{|I_1|=1,|I_2|\le 1,\az\in\{0,1,2\}}h^{I_1,I_2,\az}(t,x)\Gamma^{I_1}v\partial_\az\Gamma^{I_2}u
\ee
with the coefficients $h^{I_1,I_2,\az}(t,x)$ satisfying 
\be
|\Gamma^I(h^{I_1,I_2,\az}(t,x))|\lesssim 1\quad\mathrm{for\;any\;multi\!-\!index\;}I.
\ee
In particular, we have
\be
|P^{\gz\az\bz}\partial_\gz v\partial_\az\partial_\bz u|\lesssim\f{s^2}{t^2}|\partial_tv\partial_t\partial_tu|+\f{1}{t}\sum_{|I_1|=1,|I_2|\le 1}|\Gamma^{I_1}v||\partial\Gamma^{I_2}u|.
\ee
\item[ii)] For any sufficiently smooth functions $u,v$ supported in $\mathcal{K}$, it holds that
\be
|P^{\gz\az\bz}\partial_\gz\partial_\az v\partial_\bz u|\lesssim\f{s^2}{t^2}|\partial_t\partial_tv\partial_tu|+\f{1}{t}\sum_{|I_1|\le 1,|I_2|=1}|\partial\Gamma^{I_1}v||\Gamma^{I_2}u|.
\ee
\item[iii)] For all sufficiently smooth functions $u,v,w$ supported in $\mathcal{K}$, we have
\be
|P^{\gz\az\bz}\partial_\gz v\partial_\az u\partial_\bz w|\lesssim\f{s^2}{t^2}|\partial_tv\partial_tu\partial_tw|+\sum_{a}\l\{|\underline{\partial}_au||\partial v||\partial w|+|\partial u||\underline{\partial}_av||\partial w|+|\partial u||\partial v||\underline{\partial}_aw|\r\}.
\ee
\end{itemize}
\begin{proof}
$i)$ Denote $\xi_a=-x_a/t$. Using that $\partial_a=\xi_a\partial_t+\underline{\partial}_a$, we have
\ben
&&P^{\gz\az\bz}\partial_\gz v\partial_\az\partial_\bz u\nonumber\\
&=&P^{000}\partial_tv\partial_t\partial_tu+P^{00a}\partial_tv\partial_t\partial_au+P^{0a0}\partial_tv\partial_a\partial_tu+P^{0ab}\partial_tv\partial_a\partial_bu\nonumber\\
&+&P^{a00}\partial_av\partial_t\partial_tu+P^{a0b}\partial_av\partial_t\partial_bu+P^{ab0}\partial_av\partial_b\partial_tu+P^{abc}\partial_av\partial_b\partial_cu\nonumber\\
&=&P^{000}\partial_tv\partial_t\partial_tu+(P^{00a}+P^{0a0})\partial_tv\big(\xi_a\partial_t\partial_tu+\underline{\partial}_a\partial_tu\big)+P^{0ab}\partial_tv\big(\xi_a\partial_t+\underline{\partial}_a\big)\big(\xi_b\partial_tu+\underline{\partial}_bu\big)\nonumber\\
&+&P^{a00}\big(\xi_a\partial_tv+\underline{\partial}_av\big)\partial_t\partial_tu+(P^{a0b}+P^{ab0})\big(\xi_a\partial_tv+\underline{\partial}_av\big)\big(\xi_b\partial_t\partial_tu+\underline{\partial}_b\partial_tu\big)\nonumber\\
&+&P^{abc}\big(\xi_a\partial_tv+\underline{\partial}_av\big)\big(\xi_b\partial_t+\underline{\partial}_b\big)\big(\xi_c\partial_tu+\underline{\partial}_cu\big)\nonumber\\
&=&g(t,x)\partial_tv\partial_t\partial_tu+B(v,u),
\een
where
\ben
g(t,x):&=&P^{000}+P^{00a}\xi_a+P^{0a0}\xi_a+P^{0ab}\xi_a\xi_b+P^{a00}\xi_a+P^{a0b}\xi_a\xi_b+P^{ab0}\xi_a\xi_b+P^{abc}\xi_a\xi_b\xi_c,\nonumber\\
B(v,u):&=&(P^{00a}+P^{0a0})\partial_tv\underline{\partial}_a\partial_tu+P^{0ab}\xi_a\partial_t(\xi_b)\partial_tv\partial_tu+P^{0ab}\xi_a\partial_tv\partial_t\underline{\partial}_bu+P^{0ab}\partial_tv\underline{\partial}_a\partial_bu\nonumber\\
&+&P^{a00}\underline{\partial}_av\partial_t\partial_tu+(P^{a0b}+P^{ab0})\xi_a\partial_tv\underline{\partial}_b\partial_tu+(P^{a0b}+P^{ab0})\underline{\partial}_av\partial_b\partial_tu\nonumber\\
&+&P^{abc}\xi_a\xi_b\partial_t(\xi_c)\partial_tv\partial_tu+P^{abc}\xi_a\xi_b\partial_tv\partial_t\underline{\partial}_cu+P^{abc}\xi_a\partial_tv\underline{\partial}_b\partial_cu+P^{abc}\underline{\partial}_av\partial_b\partial_cu.
\een
We observe that $B(v,u)$ can be written as
\be
B(v,u)=\f{1}{t}\sum_{|I_1|=1,|I_2|\le 1,\az\in\{0,1,2\}}h^{I_1,I_2,\az}(t,x)\Gamma^{I_1}v\partial_\az\Gamma^{I_2}u
\ee
with the coefficients $h^{I_1,I_2,\az}(t,x)$ satisfying 
\be
|\Gamma^I(h^{I_1,I_2,\az}(t,x))|\lesssim 1\quad\mathrm{for\;any\;multi\!-\!index\;}I.
\ee
Here we use that
\be
[\partial_\az,L_b]=\dz_{\az 0}\partial_b+\dz_{\az b}\partial_t,\quad\quad[\partial_t,\underline{\partial}_a]=-\f{x_a}{t^2}\partial_t.
\ee
By the definition of $g(t,x)$ above, we also have
\be
|\Gamma^I(g(t,x))|\lesssim 1\quad\mathrm{for\;any\;multi\!-\!index\;}I.
\ee
Since $P^{\gz\az\bz}$ satisfies the null condition, we have
\be
P^{000}\f{r^3}{t^3}+P^{00a}\f{r^2}{t^2}\xi_a+P^{0a0}\f{r^2}{t^2}\xi_a+P^{0ab}\f{r}{t}\xi_a\xi_b+P^{a00}\f{r^2}{t^2}\xi_a+P^{a0b}\f{r}{t}\xi_a\xi_b+P^{ab0}\f{r}{t}\xi_a\xi_b+P^{abc}\xi_a\xi_b\xi_c=0
\ee
and therefore
\ben
g(t,x)
&=&P^{000}\l(1-\f{r^3}{t^3}\r)+P^{00a}\l(1-\f{r^2}{t^2}\r)\xi_a+P^{0a0}\l(1-\f{r^2}{t^2}\r)\xi_a+P^{0ab}\l(1-\f{r}{t}\r)\xi_a\xi_b+P^{a00}\l(1-\f{r^2}{t^2}\r)\xi_a\nonumber\\
&+&P^{a0b}\l(1-\f{r}{t}\r)\xi_a\xi_b+P^{ab0}\l(1-\f{r}{t}\r)\xi_a\xi_b.
\een
Noting that
\be
\l|1-\f{r^3}{t^3}\r|+\l|1-\f{r^2}{t^2}\r|+\l|1-\f{r}{t}\r|\lesssim\l|1-\f{r}{t}\r|\lesssim\f{s^2}{t^2},
\ee
we also have
\be
|g(t,x)|\lesssim\f{s^2}{t^2}.
\ee

$ii)$ Let $g(t,x)$ be as above. Similar to $i)$, we write 
\ben
&&P^{\gz\az\bz}\partial_\gz\partial_\az v\partial_\bz u\nonumber\\
&=&P^{000}\partial_t\partial_tv\partial_tu+P^{00a}\partial_t\partial_tv\partial_au+P^{0a0}\partial_t\partial_av\partial_tu+P^{0ab}\partial_t\partial_av\partial_bu\nonumber\\
&+&P^{a00}\partial_a\partial_tv\partial_tu+P^{a0b}\partial_a\partial_tv\partial_bu+P^{ab0}\partial_a\partial_bv\partial_tu+P^{abc}\partial_a\partial_bv\partial_cu\nonumber\\
&=&P^{000}\partial_t\partial_tv\partial_tu+P^{00a}\partial_t\partial_tv\big(\xi_a\partial_tu+\underline{\partial}_au\big)+(P^{0a0}+P^{a00})\big(\xi_a\partial_t\partial_tv+\underline{\partial}_a\partial_tv\big)\partial_tu\nonumber\\
&+&(P^{0ab}+P^{a0b})\big(\xi_a\partial_t\partial_tv+\underline{\partial}_a\partial_tv\big)\big(\xi_b\partial_tu+\underline{\partial}_bu\big)+P^{ab0}\big(\xi_a\partial_t+\underline{\partial}_a\big)\big(\xi_b\partial_tv+\underline{\partial}_bv\big)\partial_tu\nonumber\\
&+&P^{abc}\big(\xi_a\partial_t+\underline{\partial}_a\big)\big(\xi_b\partial_tv+\underline{\partial}_bv\big)\big(\xi_c\partial_tu+\underline{\partial}_cu\big)\nonumber\\
&=&g(t,x)\partial_t\partial_tv\partial_tu\nonumber\\
&+&P^{00a}\partial_t\partial_tv\underline{\partial}_au+(P^{0a0}+P^{a00})\underline{\partial}_a\partial_tv\partial_tu+(P^{0ab}+P^{a0b})\xi_a\partial_t\partial_tv\underline{\partial}_bu+(P^{0ab}+P^{a0b})\underline{\partial}_a\partial_tv\partial_bu\nonumber\\
&+&P^{ab0}\xi_a\partial_t(\xi_b)\partial_tv\partial_tu+P^{ab0}\xi_a\partial_t\underline{\partial}_bv\partial_tu+P^{ab0}\underline{\partial}_a\partial_bv\partial_tu\nonumber\\
&+&P^{abc}\xi_a\partial_t(\xi_b)\partial_tv\partial_cu+P^{abc}\xi_a\xi_b\partial_t\partial_tv\underline{\partial}_cu+P^{abc}\xi_a\partial_t\underline{\partial}_bv\partial_cu+P^{abc}\underline{\partial}_a\partial_bv\partial_cu,
\een
which implies
\be
|P^{\gz\az\bz}\partial_\gz\partial_\az v\partial_\bz u|\lesssim\f{s^2}{t^2}|\partial_t\partial_tv\partial_tu|+\f{1}{t}\sum_{|I_1|\le 1,|I_2|=1}|\partial\Gamma^{I_1}v||\Gamma^{I_2}u|.
\ee

$iii)$ Let $g(t,x)$ be as above. We have
\ben
&&P^{\gz\az\bz}\partial_\gz v\partial_\az u\partial_\bz w\nonumber\\
&=&P^{000}\partial_tv\partial_tu\partial_tw+P^{00a}\partial_tv\partial_tu\big(\xi_a\partial_tw+\underline{\partial}_aw\big)+P^{0a0}\partial_tv\big(\xi_a\partial_tu+\underline{\partial}_au\big)\partial_tw\nonumber\\
&+&P^{0ab}\partial_tv\big(\xi_a\partial_tu+\underline{\partial}_au\big)\big(\xi_b\partial_tw+\underline{\partial}_bw\big)+P^{a00}\big(\xi_a\partial_tv+\underline{\partial}_av\big)\partial_tu\partial_tw\nonumber\\
&+&P^{a0b}\big(\xi_a\partial_tv+\underline{\partial}_av\big)\partial_tu\big(\xi_b\partial_tw+\underline{\partial}_bw\big)+P^{ab0}\big(\xi_a\partial_tv+\underline{\partial}_av\big)\big(\xi_b\partial_tu+\underline{\partial}_bu\big)\partial_tw\nonumber\\
&+&P^{abc}\big(\xi_a\partial_tv+\underline{\partial}_av\big)\big(\xi_b\partial_tu+\underline{\partial}_bu\big)\big(\xi_c\partial_tw+\underline{\partial}_cw\big)\nonumber\\
&=&g(t,x)\partial_tv\partial_tu\partial_tw\nonumber\\
&+&P^{00a}\partial_tv\partial_tu\underline{\partial}_aw+P^{0a0}\partial_tv\underline{\partial}_au\partial_tw+P^{0ab}\xi_a\partial_tv\partial_tu\underline{\partial}_bw+P^{0ab}\partial_tv\underline{\partial}_au\partial_bw\nonumber\\
&+&P^{a00}\underline{\partial}_av\partial_tu\partial_tw+P^{a0b}\xi_a\partial_tv\partial_tu\underline{\partial}_bw+P^{a0b}\underline{\partial}_av\partial_tu\partial_bw+P^{ab0}\xi_a\partial_tv\underline{\partial}_bu\partial_tw\nonumber\\
&+&P^{ab0}\underline{\partial}_av\partial_bu\partial_tw+P^{abc}\xi_a\xi_b\partial_tv\partial_tu\underline{\partial}_cw+P^{abc}\xi_a\partial_tv\underline{\partial}_bu\partial_cw+P^{abc}\underline{\partial}_av\partial_bu\partial_cw,
\een
which yields
\be
|P^{\gz\az\bz}\partial_\gz v\partial_\az u\partial_\bz w|\lesssim\f{s^2}{t^2}|\partial_tv\partial_tu\partial_tw|+\sum_{a}\l\{|\underline{\partial}_au||\partial v||\partial w|+|\partial u||\underline{\partial}_av||\partial w|+|\partial u||\partial v||\underline{\partial}_aw|\r\}.
\ee
The proof is completed.
\end{proof}
\end{lem}

The Lemma below is concerned with vector fields acting on null forms; see \cite[Lemma 6.6.5]{Ho}.

\begin{lem}(See \cite{Ho})\label{lemHo}
Suppose $P^{\gz\az\bz}$ satisfies the null condition. Given any sufficiently smooth functions $v,u$, we denote $N(v,u)=P^{\gz\az\bz}\partial_\gz v\partial_\az\partial_\bz u$. Then we have 
\be
\Gamma^I(N(v,u))=\sum_{I_1+I_2\le I}N_{I;I_1,I_2}(\Gamma^{I_1}v,\Gamma^{I_2}u),
\ee
where for each $(I_1,I_2)$ with $I_1+I_2\le I$ and any sufficiently smooth functions $w,z$,
$$N_{I;I_1,I_2}(w,z)=P_{\!I;I_1,I_2}^{\gz\az\bz}\partial_\gz w\partial_\az\partial_\bz z$$ 
with $P_{\!I;I_1,I_2}^{\gz\az\bz}$ satisfying the null condition. In addition, $N_{I;I_1,I_2}(w,z)=N(w,z)$ (i.e., $P_{\!I;I_1,I_2}^{\gz\az\bz}=P^{\gz\az\bz}$) when $I_1+I_2=I$.
\end{lem}

\section{Main lemmas}\label{sE}

In this section, we present the key lemmas to be used in proving Theorem \ref{thm1}. Specifically, we establish an equality for the energy of the solution to \eqref{qEquv} up to the top order, and then introduce some nonlinear transformations for estimates of lower order energy.

\subsection{Top order energy equality}

We first show an equality for the energy of the solution to \eqref{qEquv} up to the top order.

\begin{lem}\label{lemtopE}
Let $(u,v)$ solve \eqref{qEquv} and $I\in\mathbb{N}^5$ be a multi-index with $|I|\ge 1$. For any $s\in[s_0,+\infty)$, we have 
\be
\mathcal{E}_0(\Gamma^Iu,s)+\mathcal{E}_1(\Gamma^Iv,s)
+\int_{\mathcal{H}_s}H_I\ {\rm{d}}x=\mathcal{E}_0(\Gamma^Iu,s_0)+\mathcal{E}_1(\Gamma^Iv,s_0)+\int_{\mathcal{H}_{s_0}}H_I\ {\rm{d}}x+\int_{s_0}^s\int_{\mathcal{H}_\tau}2(\tau/t)F_I\ {\rm{d}}x{\rm{d}}\tau.
\ee
Here $F_I$ is defined as
\beq\label{ExpF}
\begin{split}
&F_I:=F_{I,1}+F_{I,2}+F_{I,3},\quad\quad F_{I,1}:=F_u^{I,low}\partial_t\Gamma^Iu+F_v^{I,low}\partial_t\Gamma^Iv,\\
&F_u^{I,low}:=\sum_{I_1+I_2\le I,|I_2|\le |I|-1}\big(N_{1,I;I_1,I_2}(\Gamma^{I_1}v,\Gamma^{I_2}v)+N_{2,I;I_1,I_2}(\Gamma^{I_1}u,\Gamma^{I_2}v)\big),\\
&F_v^{I,low}:=\sum_{I_1+I_2\le I,|I_2|\le |I|-1}\big(N_{1,I;I_1,I_2}(\Gamma^{I_1}v,\Gamma^{I_2}u)+N_{2,I;I_1,I_2}(\Gamma^{I_1}u,\Gamma^{I_2}u)\big),\\
&F_{I,2}:=-\big(P_1^{\gz\az\bz}\partial_\az\partial_\gz v+P_2^{\gz\az\bz}\partial_\az\partial_\gz u\big)\big(\partial_\bz \Gamma^Iv\partial_t\Gamma^Iu+\partial_\bz \Gamma^Iu\partial_t\Gamma^Iv\big),\\
&F_{I,3}:=\f{1}{2}\big(P_1^{\gz\az\bz}\partial_t\partial_\gz v+P_2^{\gz\az\bz}\partial_t\partial_\gz u\big)\big(\partial_\bz \Gamma^Iv\partial_\az\Gamma^Iu+\partial_\bz \Gamma^Iu\partial_\az\Gamma^Iv\big),
\end{split}
\eeq
where for $i=1,2$ and any sufficiently smooth functions $w,z$, $N_{i,I;I_1,I_2}(w,z)=P_{\!i,I;I_1,I_2}^{\gz\az\bz}\partial_\gz w\partial_\az\partial_\bz z$ with $P_{\!i,I;I_1,I_2}^{\gz\az\bz}$ satisfying the null condition. The function $H_I$ is defined as
\beq\label{ExpH}
\begin{split}
H_I:=&\ \big(P_1^{\gz\az\bz}\partial_\gz v+P_2^{\gz\az\bz}\partial_\gz u\big)\big(\partial_\bz \Gamma^Iv\partial_\az\Gamma^Iu+\partial_\bz \Gamma^Iu\partial_\az\Gamma^Iv\big)\\
-&\ 2\big(P_1^{\gz\az\bz}\partial_\gz v+P_2^{\gz\az\bz}\partial_\gz u\big){\bf{n}}_\az\big(\partial_\bz \Gamma^Iv\partial_t\Gamma^Iu+\partial_\bz \Gamma^Iu\partial_t\Gamma^Iv\big)
\end{split}
\eeq
with ${\bf{n}}=(1,-x_1/t,-x_2/t)$. In addition, the following estimates hold:
\be
\begin{split}
|F_I|\lesssim&\ \f{s^2}{t^2}\sum_{\substack{|I_1|,|I_2|\le |I|\\|I_1|+|I_2|\le|I|+1}}\big(|\partial\Gamma^{I_1}v|+|\partial\Gamma^{I_1}u|\big)\big(|\partial\Gamma^{I_2}v|+|\partial\Gamma^{I_2}u|\big)\big(|\partial\Gamma^Iu|+|\partial\Gamma^Iv|\big)\\
+&\ \f{1}{t}\sum_{\substack{|I_1|,|I_2|\le |I|,|J|=1\\|I_1|+|I_2|\le|I|+1}}\big(|\Gamma^{J}\Gamma^{I_1}v|+|\Gamma^{J}\Gamma^{I_1}u|\big)\big(|\partial\Gamma^{I_2}v|+|\partial\Gamma^{I_2}u|\big)\big(|\partial\Gamma^Iu|+|\partial\Gamma^Iv|\big),\\
|H_I|\lesssim&\ \l(|\partial v|+\f{s^2}{t^2}|\partial u|\r)|\partial\Gamma^Iv||\partial\Gamma^Iu|\\
+&\ \sum_{a}\l\{|\underline{\partial}_au||\partial\Gamma^Iu||\partial\Gamma^Iv|+|\partial u|\big(|\underline{\partial}_a\Gamma^Iu|+|\underline{\partial}_a\Gamma^Iv|\big)\big(|\partial\Gamma^Iv|+|\partial\Gamma^Iu|\big)\r\}.
\end{split}
\ee
\begin{proof}
Acting the vector field $\Gamma^I$ on both sides of \eqref{qEquv} and applying Lemma \ref{lemHo}, we obtain
\beq\label{EqGzu}
-\Box \Gamma^Iu=\Gamma^IF_u=\sum_{I_1+I_2\le I}\l\{N_{1,I;I_1,I_2}(\Gamma^{I_1}v,\Gamma^{I_2}v)+N_{2,I;I_1,I_2}(\Gamma^{I_1}u,\Gamma^{I_2}v)\r\},
\eeq
\beq\label{EqGzv}
-\Box \Gamma^Iv+\Gamma^Iv=\Gamma^IF_v=\sum_{I_1+I_2\le I}\l\{N_{1,I;I_1,I_2}(\Gamma^{I_1}v,\Gamma^{I_2}u)+N_{2,I;I_1,I_2}(\Gamma^{I_1}u,\Gamma^{I_2}u)\r\},
\eeq
where for $i=1,2$ and any sufficiently smooth functions $w,z$, $N_{i,I;I_1,I_2}(w,z)=P_{\!i,I;I_1,I_2}^{\gz\az\bz}\partial_\gz w\partial_\az\partial_\bz z$ with the constant coefficients $P_{\!i,I;I_1,I_2}^{\gz\az\bz}$ satisfying the null condition, and $N_{i,I;I_1,I_2}(w,z)=N_i(w,z)$ when $I_1+I_2=I$. Multiplying \eqref{EqGzu} and \eqref{EqGzv} with $\partial_t\Gamma^Iu$ and $\partial_t\Gamma^Iv$ respectively, we have
\beq\label{Equdtu}
\f{1}{2}\partial_t(|\partial\Gamma^Iu|^2)-\partial_a(\partial^a\Gamma^Iu\partial_t\Gamma^Iu)=F_u^{I,low}\partial_t\Gamma^Iu+\big(N_1(v,\Gamma^Iv)+N_2(u,\Gamma^Iv)\big)\partial_t\Gamma^Iu,
\eeq
\beq\label{Eqvdtv}
\f{1}{2}\partial_t(|\partial\Gamma^Iv|^2+|\Gamma^Iv|^2)-\partial_a(\partial^a\Gamma^Iv\partial_t\Gamma^Iv)=F_v^{I,low}\partial_t\Gamma^Iv+\big(N_1(v,\Gamma^Iu)+N_2(u,\Gamma^Iu)\big)\partial_t\Gamma^Iv,
\eeq
where $F_u^{I,low}$ and $F_v^{I,low}$ are as in \eqref{ExpF} and we recall \eqref{dubdu}. Adding both sides of \eqref{Equdtu} and \eqref{Eqvdtv} gives
\beq\label{Equdu+vdv}
\begin{split}
&\ \f{1}{2}\partial_t(|\partial\Gamma^Iu|^2+|\partial\Gamma^Iv|^2+|\Gamma^Iv|^2)-\partial_a(\partial^a\Gamma^Iu\partial_t\Gamma^Iu+\partial^a\Gamma^Iv\partial_t\Gamma^Iv)\\
=&\ F_u^{I,low}\partial_t\Gamma^Iu+F_v^{I,low}\partial_t\Gamma^Iv+A_1+A_2,
\end{split}
\eeq
where
\ben 
A_1:&=&P_1^{\gz\az\bz}\partial_\gz v\big(\partial_\az\partial_\bz \Gamma^Iv\partial_t\Gamma^Iu+\partial_\az\partial_\bz \Gamma^Iu\partial_t\Gamma^Iv\big),\\
A_2:&=&P_2^{\gz\az\bz}\partial_\gz u\big(\partial_\az\partial_\bz \Gamma^Iv\partial_t\Gamma^Iu+\partial_\az\partial_\bz \Gamma^Iu\partial_t\Gamma^Iv\big).
\een
We have 
\beq\label{A_1}
A_1=P_1^{\gz\az\bz}\partial_\az\big(\partial_\gz v\partial_\bz \Gamma^Iv\partial_t\Gamma^Iu+\partial_\gz v\partial_\bz \Gamma^Iu\partial_t\Gamma^Iv\big)-P_1^{\gz\az\bz}\partial_\az\partial_\gz v\big(\partial_\bz \Gamma^Iv\partial_t\Gamma^Iu+\partial_\bz \Gamma^Iu\partial_t\Gamma^Iv\big)-B_1.
\eeq
Here 
\beq\label{B_1}
\begin{split}
\!\!\!\!\!\!B_1:=&\ P_1^{\gz\az\bz}\partial_\gz v\big(\partial_\bz \Gamma^Iv\partial_\az\partial_t\Gamma^Iu+\partial_\bz \Gamma^Iu\partial_\az\partial_t\Gamma^Iv\big)\\
\!\!\!\!\!\!=&\ P_1^{\gz\az\bz}\partial_t\big(\partial_\gz v\partial_\bz \Gamma^Iv\partial_\az\Gamma^Iu+\partial_\gz v\partial_\bz \Gamma^Iu\partial_\az\Gamma^Iv\big)-P_1^{\gz\az\bz}\partial_t\partial_\gz v\big(\partial_\bz \Gamma^Iv\partial_\az\Gamma^Iu+\partial_\bz \Gamma^Iu\partial_\az\Gamma^Iv\big)\\
\!\!\!\!\!\!-&\ \tilde{B}_1,
\end{split}
\eeq
where 
\beq\label{tB_1}
\tilde{B}_1:=P_1^{\gz\az\bz}\partial_\gz v\big(\partial_t\partial_\bz \Gamma^Iv\partial_\az\Gamma^Iu+\partial_t\partial_\bz \Gamma^Iu\partial_\az\Gamma^Iv\big)=B_1
\eeq
(here we use that $P_1^{\gz\az\bz}=P_1^{\gz\bz\az}$). Substituting \eqref{tB_1} into \eqref{B_1} and then by \eqref{A_1}, we derive
\be\label{ExpA_1}
\begin{split}
A_1=&\ P_1^{\gz\az\bz}\partial_\az\l\{\partial_\gz v\big(\partial_\bz \Gamma^Iv\partial_t\Gamma^Iu+\partial_\bz \Gamma^Iu\partial_t\Gamma^Iv\big)\r\}-P_1^{\gz\az\bz}\partial_\az\partial_\gz v\big(\partial_\bz \Gamma^Iv\partial_t\Gamma^Iu+\partial_\bz \Gamma^Iu\partial_t\Gamma^Iv\big)\\
-&\ \f{1}{2}P_1^{\gz\az\bz}\partial_t\l\{\partial_\gz v\big(\partial_\bz \Gamma^Iv\partial_\az\Gamma^Iu+\partial_\bz \Gamma^Iu\partial_\az\Gamma^Iv\big)\r\}+\f{1}{2}P_1^{\gz\az\bz}\partial_t\partial_\gz v\big(\partial_\bz \Gamma^Iv\partial_\az\Gamma^Iu+\partial_\bz \Gamma^Iu\partial_\az\Gamma^Iv\big).
\end{split}
\ee
In the same way, we obtain
\be\label{ExpA_2}
\begin{split}
A_2=&\ P_2^{\gz\az\bz}\partial_\az\l\{\partial_\gz u\big(\partial_\bz \Gamma^Iv\partial_t\Gamma^Iu+\partial_\bz \Gamma^Iu\partial_t\Gamma^Iv\big)\r\}-P_2^{\gz\az\bz}\partial_\az\partial_\gz u\big(\partial_\bz \Gamma^Iv\partial_t\Gamma^Iu+\partial_\bz \Gamma^Iu\partial_t\Gamma^Iv\big)\\
-& \f{1}{2}P_2^{\gz\az\bz}\partial_t\l\{\partial_\gz u\big(\partial_\bz \Gamma^Iv\partial_\az\Gamma^Iu+\partial_\bz \Gamma^Iu\partial_\az\Gamma^Iv\big)\r\}+\f{1}{2}P_2^{\gz\az\bz}\partial_t\partial_\gz u\big(\partial_\bz \Gamma^Iv\partial_\az\Gamma^Iu+\partial_\bz \Gamma^Iu\partial_\az\Gamma^Iv\big).
\end{split}
\ee
Substituting the last two equalities into \eqref{Equdu+vdv} yields
\beq\label{'Equdu+vdv}
\begin{split}
&\ \f{1}{2}\partial_t(|\partial\Gamma^Iu|^2+|\partial\Gamma^Iv|^2+|\Gamma^Iv|^2)-\partial_a(\partial^a\Gamma^Iu\partial_t\Gamma^Iu+\partial^a\Gamma^Iv\partial_t\Gamma^Iv)\\
-&\ \partial_\az\l\{\big(P_1^{\gz\az\bz}\partial_\gz v+P_2^{\gz\az\bz}\partial_\gz u\big)\big(\partial_\bz \Gamma^Iv\partial_t\Gamma^Iu+\partial_\bz \Gamma^Iu\partial_t\Gamma^Iv\big)\r\}\\
+&\ \f{1}{2}\partial_t\l\{\big(P_1^{\gz\az\bz}\partial_\gz v+P_2^{\gz\az\bz}\partial_\gz u\big)\big(\partial_\bz \Gamma^Iv\partial_\az\Gamma^Iu+\partial_\bz \Gamma^Iu\partial_\az\Gamma^Iv\big)\r\}=F_I,
\end{split}
\eeq
where $F_I$ is as in \eqref{ExpF}. Integrating \eqref{'Equdu+vdv} over $\mathcal{K}_{[s_0,s]}$, we obtain
\beq\label{estopE}
\begin{split}
&\ \int_{\mathcal{H}_s}\big(|\bar{\partial}\Gamma^Iu|^2+|\bar{\partial}\Gamma^Iv|^2+|\Gamma^Iv|^2\big){\rm{d}}x+\int_{\mathcal{H}_s}H_I\ {\rm{d}}x\\
=&\ \int_{\mathcal{H}_{s_0}}\big(|\bar{\partial}\Gamma^Iu|^2+|\bar{\partial}\Gamma^Iv|^2+|\Gamma^Iv|^2\big){\rm{d}}x+\int_{\mathcal{H}_{s_0}}H_I\ {\rm{d}}x+\int_{s_0}^s\int_{\mathcal{H}_\tau}2(\tau/t)F_I\ {\rm{d}}x{\rm{d}}\tau,
\end{split}
\eeq
where $H_I$ is as in \eqref{ExpH} and we recall \eqref{dubdu} and \eqref{DefEner}.

Recall the definition of $F_I$ in \eqref{ExpF}. For simplicity, we denote by $\partial^j$ any of the derivatives $\partial^J, |J|=j$ (recall \eqref{defdI}). By Lemma \ref{lemnf}, we have 
\ben
\!\!\!|F_{I,1}|&\lesssim&\!\!\!\sum_{\substack{|I_1|+|I_2|\le|I|\\|I_2|\le|I|-1}}\Bigg(\f{s^2}{t^2}\big(|\partial\Gamma^{I_1}v|+|\partial\Gamma^{I_1}u|\big)|\partial^2\Gamma^{I_2}v|+\f{1}{t}\sum_{\substack{|J_1|=1\\|J_2|\le 1}}\big(|\Gamma^{J_1}\Gamma^{I_1}v|+|\Gamma^{J_1}\Gamma^{I_1}u|\big)|\partial\Gamma^{J_2}\Gamma^{I_2}v|\Bigg)|\partial\Gamma^Iu|\nonumber\\
\!\!\!&+&\!\!\!\sum_{\substack{|I_1|+|I_2|\le|I|\\|I_2|\le|I|-1}}\Bigg(\f{s^2}{t^2}\big(|\partial\Gamma^{I_1}v|+|\partial\Gamma^{I_1}u|\big)|\partial^2\Gamma^{I_2}u|+\f{1}{t}\sum_{\substack{|J_1|=1\\|J_2|\le 1}}\big(|\Gamma^{J_1}\Gamma^{I_1}v|+|\Gamma^{J_1}\Gamma^{I_1}u|\big)|\partial\Gamma^{J_2}\Gamma^{I_2}u|\Bigg)|\partial\Gamma^Iv|.
\een
Lemma \ref{lemnf} also implies
\ben
|F_{I,2}|\!\!\!&\lesssim&\!\!\!\f{s^2}{t^2}\big(|\partial^2v|+|\partial^2u|\big)|\partial\Gamma^Iv||\partial\Gamma^Iu|+\f{1}{t}\sum_{|J_2|\le 1,|J_1|=1}\big(|\partial\Gamma^{J_2}v|+|\partial\Gamma^{J_2}u|\big)\big(|\Gamma^{J_1}\Gamma^Iv||\partial\Gamma^Iu|+|\Gamma^{J_1}\Gamma^Iu||\partial\Gamma^Iv|\big)
\een
and
\ben
|F_{I,3}|\!\!\!&\lesssim&\!\!\!\f{s^2}{t^2}\big(|\partial\partial_tv|+|\partial\partial_tu|\big)|\partial\Gamma^Iv||\partial\Gamma^Iu|+\sum_{a}\big(|\underline{\partial}_a\partial_tv|+|\underline{\partial}_a\partial_tu|\big)|\partial\Gamma^Iv||\partial\Gamma^Iu|\nonumber\\
\!\!\!&+&\!\!\!\sum_{a}\big(|\partial\partial_tv|+|\partial\partial_tu|\big)\big(|\underline{\partial}_a\Gamma^Iv||\partial\Gamma^Iu|+|\partial\Gamma^Iv||\underline{\partial}_a\Gamma^Iu|\big)\nonumber\\
\!\!\!&\lesssim&\!\!\!\f{s^2}{t^2}\big(|\partial^2v|+|\partial^2u|\big)|\partial\Gamma^Iv||\partial\Gamma^Iu|+\f{1}{t}\sum_{|J_2|\le 1,|J_1|=1}\big(|\partial\Gamma^{J_2}v|+|\partial\Gamma^{J_2}u|\big)\big(|\Gamma^{J_1}\Gamma^Iv||\partial\Gamma^Iu|+|\Gamma^{J_1}\Gamma^Iu||\partial\Gamma^Iv|\big).
\een
Combining the above estimates, we derive 
\ben
|F_I|\!\!\!&\lesssim&\f{s^2}{t^2}\Bigg\{\sum_{\substack{|I_1|+|I_2|\le|I|\\|I_2|\le|I|-1}}\big(|\partial\Gamma^{I_1}v|+|\partial\Gamma^{I_1}u|\big)\big(|\partial^2\Gamma^{I_2}v||\partial\Gamma^Iu|+|\partial^2\Gamma^{I_2}u||\partial\Gamma^Iv|\big)+\big(|\partial^2v|+|\partial^2u|\big)|\partial\Gamma^Iv||\partial\Gamma^Iu|\Bigg\}\nonumber\\
\!\!\!&+&\f{1}{t}\sum_{\substack{|I_1|+|I_2|\le|I|\\|I_2|\le|I|-1}}\sum_{|J_1|=1,|J_2|\le 1}\big(|\Gamma^{J_1}\Gamma^{I_1}v|+|\Gamma^{J_1}\Gamma^{I_1}u|\big)\big(|\partial\Gamma^{J_2}\Gamma^{I_2}v||\partial\Gamma^Iu|+|\partial\Gamma^{J_2}\Gamma^{I_2}u||\partial\Gamma^Iv|\big)\nonumber\\
\!\!\!&+&\f{1}{t}\sum_{|J_2|\le 1,|J_1|=1}\big(|\partial\Gamma^{J_2}v|+|\partial\Gamma^{J_2}u|\big)\big(|\Gamma^{J_1}\Gamma^Iv||\partial\Gamma^Iu|+|\Gamma^{J_1}\Gamma^Iu||\partial\Gamma^Iv|\big)\nonumber\\
\!\!\!&\lesssim&\f{s^2}{t^2}\sum_{\substack{|I_1|,|I_2|\le |I|\\|I_1|+|I_2|\le|I|+1}}\big(|\partial\Gamma^{I_1}v|+|\partial\Gamma^{I_1}u|\big)\big(|\partial\Gamma^{I_2}v|+|\partial\Gamma^{I_2}u|\big)\big(|\partial\Gamma^Iu|+|\partial\Gamma^Iv|\big)\nonumber\\
\!\!\!&+&\f{1}{t}\sum_{\substack{|I_1|,|I_2|\le |I|,|J|=1\\|I_1|+|I_2|\le|I|+1}}\big(|\Gamma^{J}\Gamma^{I_1}v|+|\Gamma^{J}\Gamma^{I_1}u|\big)\big(|\partial\Gamma^{I_2}v|+|\partial\Gamma^{I_2}u|\big)\big(|\partial\Gamma^Iu|+|\partial\Gamma^Iv|\big).
\een

For the estimate of $H_I$, we recall \eqref{ExpH} and note that ${\bf{n}}_\az=(s/t)\partial_\az s$. By Lemma \ref{lemnf}, we obtain
\ben
\!\!\!&&\!\!\!|P_2^{\gz\az\bz}\partial_\gz u\big(\partial_\bz \Gamma^Iv\partial_\az\Gamma^Iu+\partial_\bz \Gamma^Iu\partial_\az\Gamma^Iv\big)|+|P_2^{\gz\az\bz}{\bf{n}}_\az\partial_\gz u\big(\partial_\bz \Gamma^Iv\partial_t\Gamma^Iu+\partial_\bz \Gamma^Iu\partial_t\Gamma^Iv\big)|\nonumber\\
\!\!\!&\lesssim&\f{s^2}{t^2}|\partial u||\partial\Gamma^Iu||\partial\Gamma^Iv|+\sum_{a}\l\{|\underline{\partial}_au||\partial\Gamma^Iu||\partial\Gamma^Iv|+|\partial u||\underline{\partial}_a\Gamma^Iu||\partial\Gamma^Iv|+|\partial u||\partial\Gamma^Iu||\underline{\partial}_a\Gamma^Iv|\r\}\nonumber\\
\!\!\!&+&\f{s^2}{t^2}|\partial u|\l|\f{s}{t}\partial s\r||\partial\Gamma^Iv||\partial\Gamma^Iu|\nonumber\\
\!\!\!&+&\sum_{a}\l\{|\underline{\partial}_au|\l|\f{s}{t}\partial s\r||\partial\Gamma^Iv||\partial\Gamma^Iu|+|\partial u|\l|\f{s}{t}\partial s\r|\big(|\underline{\partial}_a\Gamma^Iv||\partial\Gamma^Iu|+|\underline{\partial}_a\Gamma^Iu||\partial\Gamma^Iv|\big)\r\}\nonumber\\
\!\!\!&\lesssim&\f{s^2}{t^2}|\partial u||\partial\Gamma^Iu||\partial\Gamma^Iv|+\sum_{a}\l\{|\underline{\partial}_au||\partial\Gamma^Iu||\partial\Gamma^Iv|+|\partial u||\underline{\partial}_a\Gamma^Iu||\partial\Gamma^Iv|+|\partial u||\partial\Gamma^Iu||\underline{\partial}_a\Gamma^Iv|\r\}.
\een
We also have
\be
|P_1^{\gz\az\bz}\partial_\gz v\big(\partial_\bz \Gamma^Iv\partial_\az\Gamma^Iu+\partial_\bz \Gamma^Iu\partial_\az\Gamma^Iv\big)|+|P_1^{\gz\az\bz}{\bf{n}}_\az\partial_\gz v\big(\partial_\bz \Gamma^Iv\partial_t\Gamma^Iu+\partial_\bz \Gamma^Iu\partial_t\Gamma^Iv\big)|\lesssim|\partial v||\partial\Gamma^Iv||\partial\Gamma^Iu|.
\ee
It follows that
\ben
|H_I|&\lesssim&\l(|\partial v|+\f{s^2}{t^2}|\partial u|\r)|\partial\Gamma^Iv||\partial\Gamma^Iu|\\
&+&\sum_{a}\l\{|\underline{\partial}_au||\partial\Gamma^Iu||\partial\Gamma^Iv|+|\partial u|\big(|\underline{\partial}_a\Gamma^Iu|+|\underline{\partial}_a\Gamma^Iv|\big)\big(|\partial\Gamma^Iv|+|\partial\Gamma^Iu|\big)\r\}.
\een
The proof is done.
\end{proof}
\end{lem}

\subsection{Nonlinear transforms}

To bound the lower order energy of the solution to \eqref{qEquv}, we introduce some nonlinear transformations as stated in the lemma below.

\begin{lem}\label{lemNLuv}
Let $(u,v)$ solve \eqref{qEquv} and $I\in\mathbb{N}^5$ be a multi-index. Denote 
$$\tilde{v}:=v-N_2(u,u),\quad\quad\tilde{u}:=u+\f{1}{4}N_1(v,v)+N_2(u,v).$$ 
Then 
\be\label{btv+tuR}
\begin{split}
&\ |(-\Box+1)\Gamma^I\tilde{v}|+|(-\Box)\Gamma^I\tilde{u}|\\
\lesssim&\ \f{1}{t}\sum_{|I_1|+|I_2|\le |I|+5}\big(|\Gamma^{I_1}v|+|\partial\Gamma^{I_1}u|\big)\big(|\partial\Gamma^{I_2}u|+|\Gamma^{I_2}v|\big)\\
+&\ \sum_{|I_1|+|I_2|+|I_3|\le |I|+5}\big(|\Gamma^{I_1}v|+|\partial\Gamma^{I_1}v|+|\partial\Gamma^{I_1}u|\big)\big(|\Gamma^{I_2}v|+|\partial\Gamma^{I_2}v|+|\partial\Gamma^{I_2}u|\big)|\Gamma^{I_3}v|\\
+&\ \sum_{\substack{|I_1|+|I_2|+|I_3|\le |I|+3\\|J|=1}}\l(\f{s^2}{t^2}|\partial\Gamma^{I_1}u|+\f{1}{t}|\Gamma^J\Gamma^{I_1}u|\r)|\partial\Gamma^{I_2}u||\partial\Gamma^{I_3}u|.
\end{split}
\ee
\begin{proof}
Let $\tilde{v}:=v-N_2(u,u)$. Then we have
\beq\label{Eqtv}
\begin{split}
-\Box\tilde{v}+\tilde{v}=&\ (-\Box v+v)-(-\Box)N_2(u,u)-N_2(u,u)\\
=&\ N_1(v,u)-N_2(F_u,u)-N_2(u,F_u)-2\big(N_2(\partial_tu,\partial_tu)-N_2(\partial_au,\partial^au)\big).
\end{split}
\eeq
Acting the vector field $\Gamma^I$ on both sides of \eqref{Eqtv} and applying Lemma \ref{lemHo}, we obtain
\be
(-\Box+1)\Gamma^I\tilde{v}=\mathcal{R}_1+\mathcal{R}_2,
\ee
where
\ben
\mathcal{R}_1:&=&\sum_{I_1+I_2\le I}\l\{N_{1,I;I_1,I_2}(\Gamma^{I_1}v,\Gamma^{I_2}u)-2\big(N_{2,I;I_1,I_2}(\Gamma^{I_1}\partial_tu,\Gamma^{I_2}\partial_tu)-N_{2,I;I_1,I_2}(\Gamma^{I_1}\partial_au,\Gamma^{I_2}\partial^au)\big)\r\},\nonumber\\
\mathcal{R}_2:&=&-\Gamma^I\big(N_2(F_u,u)+N_2(u,F_u)\big),
\een
where for $i=1,2$ and any sufficiently smooth functions $w,z$, $N_{i,I;I_1,I_2}(w,z)=P_{\!i,I;I_1,I_2}^{\gz\az\bz}\partial_\gz w\partial_\az\partial_\bz z$ with the constant coefficients $P_{\!i,I;I_1,I_2}^{\gz\az\bz}$ satisfying the null condition. Note that we can roughly write $N_2(F_u,u)\approx \partial F_u\partial\partial u\approx \partial(\partial v\partial\partial v+\partial u\partial\partial v)\partial\partial u,$ where we omit the constant coefficients of null forms and the subscripts of $\partial$. For simplicity, we denote by $\partial^j$ any of the derivatives $\partial^J, |J|=j$ (recall \eqref{defdI}). Then we have
\beq\label{bvR2}
\begin{split}
\!\!\!\!\!\!&\ |\mathcal{R}_2|\\
\!\!\!\!\!\!\lesssim&\ \!\!\!\sum_{|I_1|+|I_2|+|I_3|\le|I|}\Big\{\big(|\Gamma^{I_1}\partial^2 v|+|\Gamma^{I_1}\partial^2 u|\big)|\Gamma^{I_2}\partial^2 v||\Gamma^{I_3}\partial^2 u|+\big(|\Gamma^{I_1}\partial v|+|\Gamma^{I_1}\partial u|\big)|\Gamma^{I_2}\partial^3 v||\Gamma^{I_3}\partial^2 u|\\
\!\!\!\!\!\!&\quad+|\Gamma^{I_1}\partial u|\l[\big(|\Gamma^{I_2}\partial^3 v|+|\Gamma^{I_2}\partial^3 u|\big)|\Gamma^{I_3}\partial^2 v|+|\Gamma^{I_2}\partial^2 u||\Gamma^{I_3}\partial^3 v|+\big(|\Gamma^{I_2}\partial v|+|\Gamma^{I_2}\partial u|\big)|\Gamma^{I_3}\partial^4 v|\r]\Big\}\\
\!\!\!\!\!\!\lesssim&\ \sum_{|I_1|+|I_2|+|I_3|\le |I|+5}\l\{|\Gamma^{I_1}v||\Gamma^{I_2}v||\partial\Gamma^{I_3}u|+|\partial\Gamma^{I_1}u||\partial\Gamma^{I_2}u||\Gamma^{I_3}v|\r\}.
\end{split}
\eeq
By Lemmas \ref{lemnf} and \ref{lemddu},
\ben
\!\!\!\!\!\!&&|\mathcal{R}_1|\nonumber\\
\!\!\!\!\!\!&\lesssim&\!\!\!\sum_{\substack{|I_1|+|I_2|\le|I|\\|J_1|\le 1,|J_2|\le 1}}\l\{\f{s^2}{t^2}\big(|\partial\Gamma^{I_1}v||\partial^2\Gamma^{I_2}u|+|\partial\Gamma^{I_2}\partial u||\partial^2\Gamma^{I_1}\partial u|\big)+\f{1}{t}\big(|\Gamma^{J_1}\Gamma^{I_1}v||\partial\Gamma^{J_2}\Gamma^{I_2}u|+|\Gamma^{J_2}\Gamma^{I_2}\partial u||\partial\Gamma^{J_1}\Gamma^{I_1}\partial u|\big)\r\}\nonumber\\
\!\!\!\!\!\!&\lesssim&\sum_{|I_1|+|I_2|\le|I|, |J|\le 1}\Big\{\big(|\partial\Gamma^{I_1}v|+|\partial^2\Gamma^{I_1}\partial u|\big)\Big(\f{1}{t}|\partial\Gamma^{J}\Gamma^{I_2}u|+|(-\Box)\Gamma^{I_2}u|\Big)\Big\}\nonumber\\
\!\!\!\!\!\!&+&\f{1}{t}\sum_{|I_1|+|I_2|\le|I|+3}\big(|\Gamma^{I_1}v||\partial\Gamma^{I_2}u|+|\partial\Gamma^{I_1}u||\partial\Gamma^{I_2}u|\big)\nonumber\\
\!\!\!\!\!\!&\lesssim&\f{1}{t}\sum_{|I_1|+|I_2|\le|I|+3}\big(|\Gamma^{I_1}v||\partial\Gamma^{I_2}u|+|\partial\Gamma^{I_1}u||\partial\Gamma^{I_2}u|\big)+\sum_{|I_1|+|I_2|\le|I|}\big(|\partial\Gamma^{I_1}v|+|\partial^2\Gamma^{I_1}\partial u|\big)|(-\Box)\Gamma^{I_2}u|.
\een
We note that $(-\Box)\Gamma^Ju=\Gamma^{J}F_u$ for any multi-index $J$, and
\beq\label{gzKF_u}
|\Gamma^{J}F_u|\lesssim \sum_{|J_1|+|J_2|\le|J|}\big(|\Gamma^{J_1}\partial v|+|\Gamma^{J_1}\partial u|\big)|\Gamma^{J_2}\partial^2 v|
\lesssim\sum_{|J_1|+|J_2|\le|J|+2}\big(|\partial \Gamma^{J_1}v|+|\partial \Gamma^{J_1}u|\big)|\Gamma^{J_2}v|.
\eeq
It follows that
\beq\label{bvR1}
\begin{split}
|\mathcal{R}_1|\lesssim&\ \f{1}{t}\sum_{|I_1|+|I_2|\le|I|+3}\big(|\Gamma^{I_1}v||\partial\Gamma^{I_2}u|+|\partial\Gamma^{I_1}u||\partial\Gamma^{I_2}u|\big)\\
+&\ \sum_{|I_1|+|I_2|+|I_3|\le|I|+2}\big(|\partial\Gamma^{I_1}v|+|\partial^2\Gamma^{I_1}\partial u|\big)\big(|\partial \Gamma^{I_2}v|+|\partial \Gamma^{I_2}u|\big)|\Gamma^{I_3}v|\\
\lesssim&\ \f{1}{t}\sum_{|I_1|+|I_2|\le|I|+3}\big(|\Gamma^{I_1}v||\partial\Gamma^{I_2}u|+|\partial\Gamma^{I_1}u||\partial\Gamma^{I_2}u|\big)\\
+&\ \sum_{|I_1|+|I_2|+|I_3|\le|I|+5}\big(|\Gamma^{I_1}v||\Gamma^{I_2}v||\Gamma^{I_3}v|+|\partial\Gamma^{I_1}u||\Gamma^{I_2}v||\Gamma^{I_3}v|+|\partial\Gamma^{I_1}u||\partial\Gamma^{I_2}u||\Gamma^{I_3}v|\big).
\end{split}
\eeq
Combining \eqref{bvR1} and \eqref{bvR2}, we obtain
\beq\label{btvR}
\begin{split}
|(-\Box+1)\Gamma^I\tilde{v}|\lesssim&\  \f{1}{t}\sum_{|I_1|+|I_2|\le|I|+3}\big(|\Gamma^{I_1}v|+|\partial\Gamma^{I_1}u|\big)|\partial\Gamma^{I_2}u|\\
+&\ \sum_{|I_1|+|I_2|+|I_3|\le|I|+5}\!\!\big(|\Gamma^{I_1}v|+|\partial\Gamma^{I_1}u|\big)\big(|\Gamma^{I_2}v|+|\partial\Gamma^{I_2}u|\big)|\Gamma^{I_3}v|.
\end{split}
\eeq

We turn to the transformation of $u$. By \eqref{qEquv}, we have $v=\Box v+F_v$. Inserting this into the first equation in \eqref{qEquv} and applying Lemma \ref{lemnf}, we have
\beq\label{Boxu=N_boxv}
\begin{split}
-\Box u=&\ N_1(\Box v,v)+N_1(F_v,v)+N_2(u,v)\\
=&\ g_1(t,x)\big(-\partial_t\partial_t\partial_tv+\partial_t\partial_a\partial^av\big)\partial_t\partial_t v+B_1(\Box v,v)+N_1(F_v,v)+N_2(u,v),
\end{split}
\eeq
where, for any sufficiently smooth functions $w,z$,
\beq\label{'B_1(w,tw)}
B_1(w,z)=\f{1}{t}\sum_{|I_1|=1,|I_2|\le 1,\az\in\{0,1,2\}}h_1^{I_1,I_2,\az}(t,x)\Gamma^{I_1}w\partial_\az\Gamma^{I_2}z
\eeq
with $g_1(t,x)$ and $h_1^{I_1,I_2,\az}(t,x)$ satisfying
\beq\label{gzh}
|\Gamma^K(g_1(t,x))|+|\Gamma^K(h_1^{I_1,I_2,\az}(t,x))|\lesssim 1\quad\mathrm{for\;any\;multi\!-\!index\;}K.
\eeq
It follows that
\beq\label{dtvdtdtv}
g_1(t,x)\partial_t\partial_t\partial_tv\partial_t\partial_t v=\Box u+g_1(t,x)\partial_t\partial_a\partial^av\partial_t\partial_t v+B_1(\Box v,v)+N_1(F_v,v)+N_2(u,v).
\eeq
On the other hand, we have
\beq\label{BN_1(v,v)}
\begin{split}
&\ (-\Box)N_1(v,v)\\
=&\ N_1(-\Box v,v)+N_1(v,-\Box v)+2\big(N_1(\partial_tv,\partial_tv)-N_1(\partial_av,\partial^av)\big)\\
=&\ -2N_1(v,v)+N_1(F_v,v)+N_1(v,F_v)+2\big(N_1(\partial_tv,\partial_tv)-N_1(\partial_av,\partial^av)\big)\\
=&\ 2\Box u+2N_2(u,v)+N_1(F_v,v)+N_1(v,F_v)+2g_1(t,x)\big(\partial_t\partial_tv\partial_t\partial_t\partial_tv-\partial_t\partial_av\partial_t\partial_t\partial^av\big)\\
+&\ 2\big(B_1(\partial_tv,\partial_tv)-B_1(\partial_av,\partial^av)\big),
\end{split}
\eeq
where we use Lemma \ref{lemnf}, and $g_1(t,x)$ and $B_1(\cdot,\cdot)$ are as in \eqref{Boxu=N_boxv} and \eqref{'B_1(w,tw)}. Substituting \eqref{dtvdtdtv} into the right hand side of \eqref{BN_1(v,v)}, we arrive at
\beq\label{'BN_1(v,v)}
\begin{split}
&\ (-\Box)N_1(v,v)\\
=&\ 2\Box u+2N_2(u,v)+N_1(F_v,v)+N_1(v,F_v)\\
+&\ 2\Box u+2B_1(\Box v,v)+2N_1(F_v,v)+2N_2(u,v)+2g_1(t,x)\big(\partial_t\partial_a\partial^av\partial_t\partial_t v-\partial_t\partial_av\partial_t\partial_t\partial^av\big)\\
+&\ 2\big(B_1(\partial_tv,\partial_tv)-B_1(\partial_av,\partial^av)\big)\\
=&\ 4\Box u+4N_2(u,v)+3N_1(F_v,v)+N_1(v,F_v)+2B_1(\Box v,v)\\
-&\ 2g_1(t,x)\big(\underline{\partial}_a\partial_tv\partial_t\partial_t\partial^av-\partial_t\partial_tv\underline{\partial}_a\partial_t\partial^av\big)+2\big(B_1(\partial_tv,\partial_tv)-B_1(\partial_av,\partial^av)\big).
\end{split}
\eeq
We also have
\beq\label{BN_2(u,v)}
(-\Box)N_2(u,v)=N_2(F_u,v)-N_2(u,v)+N_2(u,F_v)+2\big(N_2(\partial_tu,\partial_tv)-N_2(\partial_au,\partial^av)\big).
\eeq
Let $\tilde{u}:=u+\f{1}{4}N_1(v,v)+N_2(u,v)$. Then \eqref{'BN_1(v,v)} and \eqref{BN_2(u,v)} imply
\beq\label{Eqtu}
\begin{split}
\!\!\!-\Box\tilde{u}=&\ \f{3}{4}N_1(F_v,v)+\f{1}{4}N_1(v,F_v)+N_2(F_u,v)+N_2(u,N_1(v,u))+N_2(u,N_2(u,u))+\f{1}{2}B_1(\Box v,v)\\
\!\!\!-&\ \f{1}{2}g_1(t,x)\big(\underline{\partial}_a\partial_tv\partial_t\partial_t\partial^av-\partial_t\partial_tv\underline{\partial}_a\partial_t\partial^av\big)+\f{1}{2}B_1(\partial_tv,\partial_tv)-\f{1}{2}B_1(\partial_av,\partial^av)\\
\!\!\!+&\ 2\big(N_2(\partial_tu,\partial_tv)-N_2(\partial_au,\partial^av)\big)
\end{split}
\eeq
Acting the vector field $\Gamma^I$ on both sides of \eqref{Eqtu} and applying Lemma \ref{lemHo}, we derive
\beq\label{bgztu}
\begin{split}
\!\!\!\!\!\!(-\Box)\Gamma^I\tilde{u}=&\ \mathcal{R}^I+ \sum_{I_1+I_2\le I}N_{2,I;I_1,I_2}(\Gamma^{I_1}u,\Gamma^{I_2}(N_2(u,u)))+\f{1}{2}\Gamma^I\big(B_1(\Box v,v)\big)\\
-&\ \f{1}{2}\Gamma^I\l(g_1(t,x)\big(\underline{\partial}_a\partial_tv\partial_t\partial_t\partial^av-\partial_t\partial_tv\underline{\partial}_a\partial_t\partial^av\big)\r)+\f{1}{2}\Gamma^I\big(B_1(\partial_tv,\partial_tv)-B_1(\partial_av,\partial^av)\big)\\
+&\ 2\sum_{I_1+I_2\le I}\big(N_{2,I;I_1,I_2}(\Gamma^{I_1}\partial_tu,\Gamma^{I_2}\partial_tv)-N_{2,I;I_1,I_2}(\Gamma^{I_1}\partial_au,\Gamma^{I_2}\partial^av)\big),
\end{split}
\eeq
where for any sufficiently smooth functions $w,z$, $N_{2,I;I_1,I_2}(w,z)=P_{\!2,I;I_1,I_2}^{\gz\az\bz}\partial_\gz w\partial_\az\partial_\bz z$ with $P_{\!2,I;I_1,I_2}^{\gz\az\bz}$ satisfying the null condition, and 
\beq\label{DefR^I}
\mathcal{R}^I:=\Gamma^I\l(\f{3}{4}N_1(F_v,v)+\f{1}{4}N_1(v,F_v)+N_2(F_u,v)+N_2(u,N_1(v,u))\r).
\eeq
Below we provide estimates of each term on the right hand side of \eqref{bgztu}. Note that we can write roughly $N_1(F_v,v)\approx\partial F_v\partial\partial v\approx\partial(\partial v\partial\partial u+\partial u\partial\partial u)\partial\partial v$, where we omit the constant coefficients of null forms and the subscripts of $\partial$. Hence, we have
\ben
|\Gamma^I\big(N_1(F_v,v)\big)|\!\!&\lesssim&\!\!\sum_{|I_1|+|I_2|+|I_3|\le|I|}|\Gamma^{I_1}\partial^2v|\l\{\big(|\Gamma^{I_2}\partial^2v|+|\Gamma^{I_2}\partial^2u|\big)|\Gamma^{I_3}\partial^2u|+\big(|\Gamma^{I_2}\partial v|+|\Gamma^{I_2}\partial u|\big)|\Gamma^{I_3}\partial^3u|\r\}\nonumber\\
\!\!&\lesssim&\!\!\sum_{|I_1|+|I_2|+|I_3|\le|I|+5}\big(|\Gamma^{I_1}v||\partial\Gamma^{I_2}u||\Gamma^{I_3}v|+|\partial\Gamma^{I_1}u||\partial\Gamma^{I_2}u||\Gamma^{I_3}v|\big),
\een
\ben
|\Gamma^I\big(N_1(v,F_v)\big)|\!\!&\lesssim&\!\!\sum_{|I_1|+|I_2|+|I_3|\le|I|}|\Gamma^{I_1}\partial v|\Big\{\big(|\Gamma^{I_2}\partial^3 v|+|\Gamma^{I_2}\partial^3u|\big)|\Gamma^{I_3}\partial^2u|+|\Gamma^{I_2}\partial^2 v||\Gamma^{I_3}\partial^3u|\nonumber\\
&&\quad\quad\quad\quad\quad\quad\quad\quad\quad+\big(|\Gamma^{I_2}\partial v|+|\Gamma^{I_2}\partial u|\big)|\Gamma^{I_3}\partial^4 u|\Big\}\nonumber\\
\!\!&\lesssim&\!\!\sum_{|I_1|+|I_2|+|I_3|\le|I|+5}\big(|\Gamma^{I_1}v||\Gamma^{I_2}v||\partial\Gamma^{I_3}u|+|\Gamma^{I_1}v||\partial\Gamma^{I_2}u||\partial\Gamma^{I_3}u|\big),
\een
\ben
|\Gamma^I\big(N_2(F_u,v)\big)|\!\!&\lesssim&\!\!\sum_{|I_1|+|I_2|+|I_3|\le|I|}|\Gamma^{I_1}\partial^2v|\l\{\big(|\Gamma^{I_2}\partial^2 v|+|\Gamma^{I_2}\partial^2u|\big)|\Gamma^{I_3}\partial^2v|+\big(|\Gamma^{I_2}\partial v|+|\Gamma^{I_2}\partial u|\big)|\Gamma^{I_3}\partial^3v|\r\}\nonumber\\
\!\!&\lesssim&\!\!\sum_{|I_1|+|I_2|+|I_3|\le|I|+5}\big(|\partial\Gamma^{I_1}v||\Gamma^{I_2}v||\Gamma^{I_3}v|+|\partial\Gamma^{I_1}u||\Gamma^{I_2}v||\Gamma^{I_3}v|\big)
\een
and
\ben
|\Gamma^I\big(N_2(u,N_1(v,u))\big)|\!\!&\lesssim&\!\!\sum_{|I_1|+|I_2|+|I_3|\le|I|}|\Gamma^{I_1}\partial u|\big(|\Gamma^{I_2}\partial^3v||\Gamma^{I_3}\partial^2u|+|\Gamma^{I_2}\partial^2v||\Gamma^{I_3}\partial^3u|+|\Gamma^{I_2}\partial v||\Gamma^{I_3}\partial^4u|\big)\nonumber\\
\!\!&\lesssim&\!\!\sum_{|I_1|+|I_2|+|I_3|\le |I|+4}|\partial\Gamma^{I_1}u||\partial\Gamma^{I_2}u||\Gamma^{I_3}v|,
\een
which imply that
\beq\label{bu1}
\mathcal{R}^I\lesssim\ \sum_{|I_1|+|I_2|+|I_3|\le|I|+5}\big(|\Gamma^{I_1}v||\Gamma^{I_2}v||\partial\Gamma^{I_3}u|+|\partial\Gamma^{I_1}u||\partial\Gamma^{I_2}u||\Gamma^{I_3}v|+|\partial\Gamma^{I_1}v||\Gamma^{I_2}v||\Gamma^{I_3}v|\big).
\eeq
We write $N_2(u,u)$ roughly as $\partial u\partial\partial u$ (here we omit the coefficients of null forms and the subscripts of $\partial$). By Lemma \ref{lemnf}, we have
\beq\label{N_2(u,N_2u)}
\begin{split}
&\ \sum_{I_1+I_2\le I}\l|N_{2,I;I_1,I_2}\big(\Gamma^{I_1}u,\Gamma^{I_2}(N_2(u,u))\big)\r|\\
\lesssim&\ \sum_{|I_1|+|I_2|\le|I|}\Bigg\{\f{s^2}{t^2}|\partial\Gamma^{I_1}u||\partial^2\Gamma^{I_2}(\partial u\partial^2u)|+\f{1}{t}\sum_{|J_1|=1,|J_2|\le 1}|\Gamma^{J_1}\Gamma^{I_1}u||\partial\Gamma^{J_2}\Gamma^{I_2}(\partial u\partial^2u)|\Bigg\}\\
\lesssim&\ \sum_{\substack{|I_1|+|I_2|+|I_3|\le|I|+3\\|J|=1}}\l(\f{s^2}{t^2}|\partial\Gamma^{I_1}u|+\f{1}{t}|\Gamma^J\Gamma^{I_1}u|\r)|\partial\Gamma^{I_2}u||\partial\Gamma^{I_3}u|.
\end{split}
\eeq
By \eqref{'B_1(w,tw)} and \eqref{gzh}, we have
\beq\label{gzB_1bvv}
\begin{split}
|\Gamma^I\big(B_1(\Box v,v)\big)|=&\ \Bigg|\Gamma^I\Bigg(\f{1}{t}\sum_{|J_1|=1,|J_2|\le 1,\az\in\{0,1,2\}}h_1^{J_1,J_2,\az}(t,x)\Gamma^{J_1}\Box v\partial_\az\Gamma^{J_2}v\Bigg)\Bigg|\\
\lesssim&\ \f{1}{t}\sum_{|I_1|+|I_2|\le|I|,|J_1|,|J_2|\le 1}|\Gamma^{I_1}\Gamma^{J_1}\partial^2v||\Gamma^{I_2}\partial\Gamma^{J_2}v|\lesssim\f{1}{t}\sum_{|I_1|+|I_2|\le|I|+5}|\Gamma^{I_1}v||\Gamma^{I_2}v|
\end{split}
\eeq
and similarly
\beq\label{gzB_1dtv}
\begin{split}
|\Gamma^I\big(B_1(\partial_tv,\partial_tv)\big)|+|\Gamma^I\big(B_1(\partial_av,\partial^av)\big)|
\lesssim&\ \f{1}{t}\sum_{|I_1|+|I_2|\le|I|,|J_1|,|J_2|\le 1}|\Gamma^{I_1}\Gamma^{J_1}\partial v||\Gamma^{I_2}\partial\Gamma^{J_2}\partial v|\\
\lesssim&\ \f{1}{t}\sum_{|I_1|+|I_2|\le|I|+5}|\Gamma^{I_1}v||\Gamma^{I_2}v|.
\end{split}
\eeq
Using \eqref{gzh} again, we obtain
\beq\label{gzg1bdv}
\l|\Gamma^I\l(g_1(t,x)\big(\underline{\partial}_a\partial_tv\partial_t\partial_t\partial^av-\partial_t\partial_tv\underline{\partial}_a\partial_t\partial^av\big)\r)\r|\lesssim\f{1}{t}\sum_{|I_1|+|I_2|\le|I|+5}|\Gamma^{I_1}v||\Gamma^{I_2}v|.
\eeq
By Lemmas \ref{lemnf} and \ref{lemddu}, we have
\beq\label{N_2ddudv}
\begin{split}
&\ \sum_{I_1+I_2\le I}\big(|N_{2,I;I_1,I_2}(\Gamma^{I_1}\partial_tu,\Gamma^{I_2}\partial_tv)|+|N_{2,I;I_1,I_2}(\Gamma^{I_1}\partial_au,\Gamma^{I_2}\partial^av)|\big)\\
\lesssim&\ \f{s^2}{t^2}\sum_{|I_1|+|I_2|\le|I|}|\partial\Gamma^{I_1}\partial u||\partial^2\Gamma^{I_2}\partial v|+\f{1}{t}\sum_{\substack{|I_1|+|I_2|\le|I|\\|J_1|=1,|J_2|\le 1}}|\Gamma^{J_1}\Gamma^{I_1}\partial u||\partial\Gamma^{J_2}\Gamma^{I_2}\partial v|\\
\lesssim&\ \f{s^2}{t^2}\sum_{|I_1|+|I_2|\le|I|+3}|\partial\partial\Gamma^{I_1}u||\Gamma^{I_2}v|+\f{1}{t}\sum_{|I_1|+|I_2|\le|I|+4}|\partial\Gamma^{I_1}u||\Gamma^{I_2}v|\\
\lesssim&\ \sum_{\substack{|I_1|+|I_2|\le|I|+3\\|J|\le 1}}\l(\f{1}{t}|\partial\Gamma^J\Gamma^{I_1}u|+|(-\Box)\Gamma^{I_1}u|\r)|\Gamma^{I_2}v|+\f{1}{t}\sum_{|I_1|+|I_2|\le|I|+4}|\partial\Gamma^{I_1}u||\Gamma^{I_2}v|\\
\lesssim&\ \f{1}{t}\sum_{|I_1|+|I_2|\le|I|+4}|\partial\Gamma^{I_1}u||\Gamma^{I_2}v|+\sum_{|I_1|+|I_2|+|I_3|\le|I|+5}\big(|\partial\Gamma^{I_1}v|+|\partial\Gamma^{I_1}u|\big)|\Gamma^{I_2}v||\Gamma^{I_3}v|,
\end{split}
\eeq
where we use \eqref{gzKF_u}. Combining \eqref{bgztu}-\eqref{N_2ddudv}, we derive
\be\label{btuR}
\begin{split}
\!\!\!|(-\Box)\Gamma^I\tilde{u}|\lesssim&\ \sum_{|I_1|+|I_2|+|I_3|\le|I|+5}\big(|\partial\Gamma^{I_1}v|+|\partial\Gamma^{I_1}u|\big)\big(|\Gamma^{I_2}v|+|\partial\Gamma^{I_2}u|\big)|\Gamma^{I_3}v|\\
\!\!\!+&\ \sum_{\substack{|I_1|+|I_2|+|I_3|\le|I|+3\\|J|=1}}\l(\f{s^2}{t^2}|\partial\Gamma^{I_1}u|+\f{1}{t}|\Gamma^J\Gamma^{I_1}u|\r)|\partial\Gamma^{I_2}u||\partial\Gamma^{I_3}u|+\ \f{1}{t}\sum_{|I_1|+|I_2|\le|I|+5}\big(|\Gamma^{I_1}v|+|\partial\Gamma^{I_1}u|\big)|\Gamma^{I_2}v|.
\end{split}
\ee
This together with \eqref{btvR} yields
\ben
&&|(-\Box+1)\Gamma^I\tilde{v}|+|(-\Box)\Gamma^I\tilde{u}|\nonumber\\
&\lesssim&\f{1}{t}\sum_{|I_1|+|I_2|\le |I|+5}\big(|\Gamma^{I_1}v|+|\partial\Gamma^{I_1}u|\big)\big(|\partial\Gamma^{I_2}u|+|\Gamma^{I_2}v|\big)\nonumber\\
&+&\sum_{|I_1|+|I_2|+|I_3|\le |I|+5}\big(|\Gamma^{I_1}v|+|\partial\Gamma^{I_1}v|+|\partial\Gamma^{I_1}u|\big)\big(|\Gamma^{I_2}v|+|\partial\Gamma^{I_2}v|+|\partial\Gamma^{I_2}u|\big)|\Gamma^{I_3}v|\nonumber\\
&+&\sum_{\substack{|I_1|+|I_2|+|I_3|\le |I|+3\\|J|=1}}\l(\f{s^2}{t^2}|\partial\Gamma^{I_1}u|+\f{1}{t}|\Gamma^J\Gamma^{I_1}u|\r)|\partial\Gamma^{I_2}u||\partial\Gamma^{I_3}u|.
\een
The proof is completed.
\end{proof}
\end{lem}

\section{Proof of the global existence}\label{sM}

In this section we prove the global existence result in Theorem \ref{thm1}.

\subsection{Bootstrap assumption}

${\bf{Bootstrap\ setting}.}$ Let $N\ge 14$ be an integer and $0<\dz\ll 1$. We assume the following energy bound for the solution $(u,v)$ to \eqref{qEquv}-\eqref{qini} on $[s_0,s_*)$:
\beq\label{Bsuv}
\!\!\sum_{|I|\le N}\!\big\{[\mathcal{E}_0(\Gamma^Iu,s)]^{1/2}+[\mathcal{E}_1(\Gamma^Iv,s)]^{1/2}\big\}\cdot s^{-\dz}+\sum_{|I|\le N-5}\!\!\big\{[\mathcal{E}_0(\Gamma^Iu,s)]^{1/2}+[\mathcal{E}_1(\Gamma^Iv,s)]^{1/2}\big\}\le C_1\ez
\eeq
with $C_1>1$ some large constant to be determined later, $0<\ez\ll 1$ the size of the initial data and
\beq\label{Defs_*}
s_*:=\sup\{s>s_0: \eqref{Bsuv}\;\mathrm{holds\;on\;}[s_0,s]\}.
\eeq

\begin{prop}\label{propim}
There exist some constants $C_1$ sufficiently large and $0<\ez_0\ll C_1^{-1}$ sufficiently small such that, for any $0<\ez<\ez_0$, if $(u,v)$ is a solution to \eqref{qEquv}-\eqref{qini} and satisfies \eqref{Bsuv} on the hyperbolic time interval $[s_0,s_1]$, then for $s\in[s_0,s_1]$ we have
\be
\sum_{|I|\le N}\big\{[\mathcal{E}_0(\Gamma^Iu,s)]^{1/2}+[\mathcal{E}_1(\Gamma^Iv,s)]^{1/2}\big\}\cdot s^{-\dz}+\sum_{|I|\le N-5}\big\{[\mathcal{E}_0(\Gamma^Iu,s)]^{1/2}+[\mathcal{E}_1(\Gamma^Iv,s)]^{1/2}\big\}\le \f{1}{2}C_1\ez.
\ee
\end{prop}
In the above proposition $s_1$ is arbitrary, hence $s_*=+\infty$ where $s_*$ is as in \eqref{Defs_*}, which implies that the solution $(u,v)$ exists globally in time and satisfies \eqref{Bsuv} for any $s\in[s_0,+\infty)$. Below we provide the proof of Proposition \ref{propim}. In the sequel, the implied constants in "$\lesssim$" do not depend on the constants $C_1$ and $\ez$ appearing in the bootstrap assumption \eqref{Bsuv}.

Let $(u,v)$ be a solution to \eqref{qEquv}-\eqref{qini} and satisfy \eqref{Bsuv} for $s\in[s_0,s_1]$. By \eqref{DefEner},\eqref{s/tduE} and Lemmas \ref{sC: DL2.4} and \ref{sC: DL2.3}, we have the following $L^2$-type and pointwise estimates for $s\in[s_0,s_1]$:
\beq\label{L2Linf}
\l\{\begin{split}
&\ \sum_{|I|\le N,a\in\{1,2\}}\|(s/t)\big(|\partial\Gamma^Iu|+|\partial\Gamma^Iv|\big)+\big(|\underline{\partial}_a\Gamma^Iu|+|\underline{\partial}_a\Gamma^Iv|+|\Gamma^Iv|\big)\|_{L^2_f(\mathcal{H}_s)}\lesssim C_1\ez s^\dz,\\
&\ \sum_{|I|\le N-7,a\in\{1,2\}}\sup_{\mathcal{H}_s}\l\{s\big(|\partial\Gamma^Iu|+|\partial\Gamma^Iv|\big)+t\big(|\underline{\partial}_a\Gamma^Iu|+|\underline{\partial}_a\Gamma^Iv|+|\Gamma^Iv|\big)\r\}\lesssim C_1\ez.
\end{split}\r.
\eeq
By \eqref{L2Linf}, we also obtain
\beq\label{gzL2Linf}
\sum_{\substack{|I|\le N\\|J|=1}}\|t^{-1}\big(|\Gamma^J\Gamma^Iu|+|\Gamma^J\Gamma^Iv|\big)\|_{L^2_f(\mathcal{H}_s)}\cdot s^{-\dz}+\sum_{\substack{|I|\le N-7\\|J|=1}}\|\ \!|\Gamma^J\Gamma^Iu|+|\Gamma^J\Gamma^Iv|\ \!\|_{L^\infty(\mathcal{H}_s)}\lesssim C_1\ez.
\eeq

\subsection{Improved estimates of energy up to the top order}

In this subsection we show refined estimates of energy of the solution $(u,v)$ up to the top order. 

By Lemma \ref{lemtopE}, for $|I|\le N$, we have
\beq\label{E0us+E1vs=}
\begin{split}
\!\!\!\mathcal{E}_0(\Gamma^Iu,s)+\mathcal{E}_1(\Gamma^Iv,s)+\int_{\mathcal{H}_s}H_I\ {\rm{d}}x
=&\ \ \mathcal{E}_0(\Gamma^Iu,s_0)+\mathcal{E}_1(\Gamma^Iv,s_0)+\int_{\mathcal{H}_{s_0}}H_I\ {\rm{d}}x\\
\!\!\!+&\ \ \int_{s_0}^s\int_{\mathcal{H}_\tau}2(\tau/t)F_I\ {\rm{d}}x{\rm{d}}\tau,
\end{split}
\eeq
where $F_I$, $H_I$ are as in \eqref{ExpF} and \eqref{ExpH} respectively. For any $\tau\in[s_0,s]$, Lemma \ref{lemtopE} implies that on $\mathcal{H}_\tau$
\ben
|F_I|\!\!\!&\lesssim&\!\!\!\f{\tau^2}{t^2}\sum_{\substack{|I_1|,|I_2|\le N\\|I_1|+|I_2|\le N+1}}\big(|\partial\Gamma^{I_1}v|+|\partial\Gamma^{I_1}u|\big)\big(|\partial\Gamma^{I_2}v|+|\partial\Gamma^{I_2}u|\big)\big(|\partial\Gamma^Iu|+|\partial\Gamma^Iv|\big)\nonumber\\
\!\!\!&+&\!\!\!\f{1}{t}\sum_{\substack{|I_1|,|I_2|\le N,|J|=1\\|I_1|+|I_2|\le N+1}}\big(|\Gamma^{J}\Gamma^{I_1}v|+|\Gamma^{J}\Gamma^{I_1}u|\big)\big(|\partial\Gamma^{I_2}v|+|\partial\Gamma^{I_2}u|\big)\big(|\partial\Gamma^Iu|+|\partial\Gamma^Iv|\big).
\een
Hence by \eqref{L2Linf} and \eqref{gzL2Linf}, we have
\beq\label{tau/tFL1}
\begin{split}
\!\!\!\l\|\f{\tau}{t}F_I\r\|_{L^1_f(\mathcal{H}_\tau)}\lesssim&\ \sum_{\substack{|I_1|\le N-7\\|I_2|\le N}}\Big\{\|\tau\big(|\partial\Gamma^{I_1}v|+|\partial\Gamma^{I_1}u|\big)\|_{L^\infty(\mathcal{H}_\tau)}\cdot\|(\tau/t)\big(|\partial\Gamma^{I_2}v|+|\partial\Gamma^{I_2}u|\big)\|_{L^2_f(\mathcal{H}_\tau)}\\
\!\!\!&\quad\quad\quad\; \cdot\ \|(\tau/t)\big(|\partial\Gamma^{I}v|+|\partial\Gamma^{I}u|\big)\|_{L^2_f(\mathcal{H}_\tau)}\cdot\tau^{-1}\Big\}\\
\!\!\!+&\ \sum_{\substack{|I_1|\le N-7\\|I_2|\le N,|J|=1}}\Big\{\|\ \!|\Gamma^J\Gamma^{I_1}v|+|\Gamma^J\Gamma^{I_1}u|\ \!\|_{L^\infty(\mathcal{H}_\tau)}\cdot\|(\tau/t)\big(|\partial\Gamma^{I_2}v|+|\partial\Gamma^{I_2}u|\big)\|_{L^2_f(\mathcal{H}_\tau)}\\
\!\!\!&\quad\quad\quad\quad\; \cdot\ \|(\tau/t)\big(|\partial\Gamma^{I}v|+|\partial\Gamma^{I}u|\big)\|_{L^2_f(\mathcal{H}_\tau)}\cdot\tau^{-1}\Big\}\\
\!\!\!+&\ \sum_{\substack{|I_2|\le N-7\\|I_1|\le N,|J|=1}}\Big\{\|t^{-1}\big(|\Gamma^J\Gamma^{I_1}v|+|\Gamma^J\Gamma^{I_1}u|\big)\|_{L^2_f(\mathcal{H}_\tau)}\cdot\|\tau\big(|\partial\Gamma^{I_2}v|+|\partial\Gamma^{I_2}u|\big)\|_{L^\infty(\mathcal{H}_\tau)}\\
\!\!\!&\quad\quad\quad\quad\;\cdot\ \|(\tau/t)\big(|\partial\Gamma^{I}v|+|\partial\Gamma^{I}u|\big)\|_{L^2_f(\mathcal{H}_\tau)}\cdot\tau^{-1}\Big\}\\
\!\!\!\lesssim&\ (C_1\ez)^3\tau^{-1+2\dz},
\end{split}
\eeq
where we use that $N\ge 14$. We also obtain from Lemma \ref{lemtopE} that
\ben
|H_I|&\lesssim&|\partial v||\partial\Gamma^Iv||\partial\Gamma^Iu|+\f{s^2}{t^2}|\partial u||\partial\Gamma^Iu||\partial\Gamma^Iv|\nonumber\\
&+&\sum_{a}\l\{|\underline{\partial}_au||\partial\Gamma^Iu||\partial\Gamma^Iv|+|\partial u|\big(|\underline{\partial}_a\Gamma^Iu|+|\underline{\partial}_a\Gamma^Iv|\big)\big(|\partial\Gamma^Iv|+|\partial\Gamma^Iu|\big)\r\},
\een
which implies
\ben
\|H_I\|_{L^1_f(\mathcal{H}_s)}&\lesssim&\sum_{a}\|t|\partial v|+|\partial u|+t|\underline{\partial}_au|\ \!\|_{L^\infty(\mathcal{H}_s)}\cdot\|(s/t)|\partial\Gamma^Iv|\ \!\|_{L^2_f(\mathcal{H}_s)}\cdot\|(s/t)|\partial\Gamma^Iu|\ \!\|_{L^2_f(\mathcal{H}_s)}\nonumber\\
&+&\sum_{a}\|s|\partial u|\ \!\|_{L^\infty(\mathcal{H}_s)}\cdot\|\ \!|\underline{\partial}_a\Gamma^Iv|+|\underline{\partial}_a\Gamma^Iu|\ \!\|_{L^2_f(\mathcal{H}_s)}\cdot\|(s/t)\big(|\partial\Gamma^Iv|+|\partial\Gamma^Iu|\big)\|_{L^2_f(\mathcal{H}_s)}\nonumber\\
&\lesssim&C_1\ez\l\{\mathcal{E}_0(\Gamma^Iu,s)+\mathcal{E}_1(\Gamma^Iv,s)\r\}\le\f{1}{2}\l\{\mathcal{E}_0(\Gamma^Iu,s)+\mathcal{E}_1(\Gamma^Iv,s)\r\},
\een
where we choose $0<\ez\ll C_1^{-1}$ sufficiently small. Combining this with \eqref{E0us+E1vs=} and \eqref{tau/tFL1}, we conclude that
\beq\label{imTuv}
\mathcal{E}_0(\Gamma^Iu,s)+\mathcal{E}_1(\Gamma^Iv,s)\lesssim\ez^2+(C_1\ez)^3s^{2\dz},\quad|I|\le N.
\eeq
Hence we have strictly improved the estimates of $\sum_{|I|\le N}\l\{[\mathcal{E}_0(\Gamma^Iu,s)]^{1/2}+[\mathcal{E}_1(\Gamma^Iv,s)]^{1/2}\r\}$ in \eqref{Bsuv} if we choose $C_1>1$ sufficiently large and $0<\ez\ll C_1^{-1}$ sufficiently small.

\subsection{Improved estimates of lower order energy}

In this subsection we show refined estimates of lower order energy of the solution $(u,v)$. 

Let $\tilde{v}:=v-N_2(u,u)$ and $\tilde{u}:=u+\f{1}{4}N_1(v,v)+N_2(u,v)$. By Lemma \ref{lemNLuv}, we have
\beq\label{'btv+tuR}
|(-\Box+1)\Gamma^I\tilde{v}|+|(-\Box)\Gamma^I\tilde{u}|\lesssim Q_1+Q_2+Q_3,\quad|I|\le N-5,
\eeq
where
\ben
Q_1&=&\f{1}{t}\sum_{|I_1|+|I_2|\le N}\big(|\Gamma^{I_1}v|+|\partial\Gamma^{I_1}u|\big)\big(|\partial\Gamma^{I_2}u|+|\Gamma^{I_2}v|\big)\nonumber\\
Q_2&=&\sum_{|I_1|+|I_2|+|I_3|\le N}\big(|\Gamma^{I_1}v|+|\partial\Gamma^{I_1}v|+|\partial\Gamma^{I_1}u|\big)\big(|\Gamma^{I_2}v|+|\partial\Gamma^{I_2}v|+|\partial\Gamma^{I_2}u|\big)|\Gamma^{I_3}v|\nonumber\\
Q_3&=&\sum_{\substack{|I_1|+|I_2|+|I_3|\le N-2\\|J|=1}}\l(\f{s^2}{t^2}|\partial\Gamma^{I_1}u|+\f{1}{t}|\Gamma^J\Gamma^{I_1}u|\r)|\partial\Gamma^{I_2}u||\partial\Gamma^{I_3}u|.
\een
For $|I|\le N-5$, Proposition \ref{sC: DL5.1} implies
\beq\label{E0tu_E1tv}
\begin{split}
[\mathcal{E}_0(\Gamma^I\tilde{u},s)]^{1/2}+[\mathcal{E}_1(\Gamma^I\tilde{v},s)]^{1/2}\lesssim&\;\   [\mathcal{E}_0(\Gamma^I\tilde{u},s_0)]^{1/2}+[\mathcal{E}_1(\Gamma^I\tilde{v},s_0)]^{1/2}\\
+&\ \int_{s_0}^s\|Q_1+Q_2+Q_3\|_{L^2_f(\mathcal{H}_\tau)}{\rm{d}\tau}.
\end{split}
\eeq
For any $\tau\in[s_0,s]$, by \eqref{L2Linf}, we derive
\beq\label{EsQ1}
\begin{split}
\|Q_1\|_{L^2_f(\mathcal{H}_\tau)}\lesssim&\ \sum_{\substack{|I_1|\le N-7\\|I_2|\le N}}\|\tau\big(|\Gamma^{I_1}v|+|\partial\Gamma^{I_1}u|\big)\|_{L^\infty(\mathcal{H}_\tau)}\cdot\|(\tau/t)\big(|\partial\Gamma^{I_2}u|+|\Gamma^{I_2}v|\big)\|_{L^2_f(\mathcal{H}_\tau)}\cdot\tau^{-2}\\
\lesssim&\ (C_1\ez)^2\tau^{-2+\dz},
\end{split}
\eeq
\beq\label{EsQ2}
\begin{split}
\|Q_2\|_{L^2_f(\mathcal{H}_\tau)}\lesssim&\ \sum_{\substack{|I_1|,|I_2|\le N-7\\|I_3|\le N}}\Big\{\|\tau\big(|\Gamma^{I_1}v|+|\partial\Gamma^{I_1}v|+|\partial\Gamma^{I_1}u|\big)\|_{L^\infty(\mathcal{H}_\tau)}\\
&\quad\quad\quad\quad\quad\cdot\ \|\tau\big(|\Gamma^{I_2}v|+|\partial\Gamma^{I_2}v|+|\partial\Gamma^{I_2}u|\big)\|_{L^\infty(\mathcal{H}_\tau)}\cdot\|\ \!|\Gamma^{I_3}v|\ \!\|_{L^2_f(\mathcal{H}_\tau)}\cdot\tau^{-2}\Big\}\\
+&\ \sum_{\substack{|I_2|,|I_3|\le N-7\\|I_1|\le N}}\Big\{\|(\tau/t)\big(|\Gamma^{I_1}v|+|\partial\Gamma^{I_1}v|+|\partial\Gamma^{I_1}u|\big)\|_{L^2_f(\mathcal{H}_\tau)}\\
&\quad\quad\quad\quad\quad\cdot\ \|\tau\big(|\Gamma^{I_2}v|+|\partial\Gamma^{I_2}v|+|\partial\Gamma^{I_2}u|\big)\|_{L^\infty(\mathcal{H}_\tau)}\cdot\|t|\Gamma^{I_3}v|\ \!\|_{L^\infty(\mathcal{H}_\tau)}\cdot\tau^{-2}\Big\}\\
\lesssim&\ (C_1\ez)^3\tau^{-2+\dz}.
\end{split}
\eeq
By \eqref{L2Linf} and \eqref{gzL2Linf}, we also obtain
\beq\label{EsQ3}
\begin{split}
\!\!\!\|Q_3\|_{L^2_f(\mathcal{H}_\tau)}\lesssim&\ \sum_{\substack{|I_2|,|I_3|\le N-7\\|I_1|\le N,|J|=1}}\Big\{\|(\tau/t)|\partial\Gamma^{I_1}u|+t^{-1}|\Gamma^J\Gamma^{I_1}u|\ \!\|_{L^2_f(\mathcal{H}_\tau)}\cdot\|\tau|\partial\Gamma^{I_2}u|\ \!\|_{L^\infty(\mathcal{H}_\tau)}\\
\!\!\!&\quad\quad\quad\quad\quad\cdot\|\tau|\partial\Gamma^{I_3}u|\ \!\|_{L^\infty(\mathcal{H}_\tau)}\cdot\tau^{-2}\Big\}\\
\!\!\!+&\ \sum_{\substack{|I_1|,|I_2|\le N-7\\|I_3|\le N,|J|=1}}\|\ \!|\Gamma^J\Gamma^{I_1}u|\ \!\|_{L^\infty(\mathcal{H}_\tau)}\cdot\|\tau|\partial\Gamma^{I_2}u|\ \!\|_{L^\infty(\mathcal{H}_\tau)}\cdot\|(\tau/t)|\partial\Gamma^{I_3}u|\ \!\|_{L^2_f(\mathcal{H}_\tau)}\cdot\tau^{-2}\\
\!\!\!\lesssim&\ (C_1\ez)^3\tau^{-2+\dz}.
\end{split}
\eeq
Combining \eqref{E0tu_E1tv}-\eqref{EsQ3}, we obtain
\beq\label{imtutv}
[\mathcal{E}_0(\Gamma^I\tilde{u},s)]^{1/2}+[\mathcal{E}_1(\Gamma^I\tilde{v},s)]^{1/2}\lesssim\ez+(C_1\ez)^2,\quad |I|\le N-5.
\eeq
Note that we can write roughly $N_1(v,v)\approx\partial v\partial\partial v$, $N_2(u,v)\approx\partial u\partial\partial v$ and $N_2(u,u)\approx\partial u\partial\partial u$ (here we omit the constant coefficients of null forms and the subscripts of $\partial$). Hence for $|I|\le N-5$, we obtain
\beq\label{E_0u-tu}
\begin{split}
&\ [\mathcal{E}_0(\Gamma^I(u-\tilde{u}),s)]^{1/2}\lesssim\|\ \!|\partial\Gamma^I(N_1(v,v))|+|\partial\Gamma^I(N_2(u,v))|\ \!\|_{L^2_f(\mathcal{H}_s)}\\
\lesssim&\ \sum_{|I_1|+|I_2|\le|I|}\|\big(|\partial\Gamma^{I_1}\partial v|+|\partial\Gamma^{I_1}\partial u|\big)|\Gamma^{I_2}\partial^2v|+\big(|\Gamma^{I_1}\partial v|+|\Gamma^{I_1}\partial u|\big)|\partial\Gamma^{I_2}\partial^2v|\ \!\|_{L^2_f(\mathcal{H}_s)}\\
\lesssim&\ \sum_{|I_1|+|I_2|\le|I|+3}\|\big(|\partial\Gamma^{I_1}v|+|\partial\Gamma^{I_1}u|\big)|\Gamma^{I_2}v|\ \!\|_{L^2_f(\mathcal{H}_s)}\\
\lesssim&\ \sum_{|I_1|\le N-7,|I_2|\le N}\|s\big(|\partial\Gamma^{I_1}v|+|\partial\Gamma^{I_1}u|\big)\|_{L^\infty(\mathcal{H}_s)}\cdot\|\ \!|\Gamma^{I_2}v|\ \!\|_{L^2_f(\mathcal{H}_s)}\cdot s^{-1}\\
\lesssim&\ \sum_{|I_2|\le N-7,|I_1|\le N}\|(s/t)\big(|\partial\Gamma^{I_1}v|+|\partial\Gamma^{I_1}u|\big)\|_{L^2_f(\mathcal{H}_s)}\cdot\|t|\Gamma^{I_2}v|\ \!\|_{L^\infty(\mathcal{H}_s)}\cdot s^{-1}\\
\lesssim&\  (C_1\ez)^2s^{-1+\dz}\lesssim(C_1\ez)^2,
\end{split}
\eeq
where we use \eqref{L2Linf}. For $|I|\le N-5$, we also have (recall \eqref{DefEner})
\be
[\mathcal{E}_1(\Gamma^I(v-\tilde{v}),s)]^{1/2}\lesssim\sum_{a}\|(s/t)|\partial\Gamma^I(N_2(u,u))|+|\underline{\partial}_a\Gamma^I(N_2(u,u))|+|\Gamma^I(N_2(u,u))|\ \!\|_{L^2_f(\mathcal{H}_s)}.
\ee
We compute 
\ben
\!\!\!&&\!\!\!\sum_{a}\|(s/t)|\partial\Gamma^I(N_2(u,u))|+|\underline{\partial}_a\Gamma^I(N_2(u,u))|\ \!\|_{L^2_f(\mathcal{H}_s)}\nonumber\\
\!\!\!&\lesssim&\!\!\!\sum_{\substack{|I_1|+|I_2|\le|I|\\a\in\{1,2\}}}\!\!\!\|(s/t)\big(|\partial\Gamma^{I_1}\partial u||\Gamma^{I_2}\partial^2u|+|\Gamma^{I_1}\partial u||\partial\Gamma^{I_2}\partial^2u|\big)+|\underline{\partial}_a\Gamma^{I_1}\partial u||\Gamma^{I_2}\partial^2u|+|\Gamma^{I_1}\partial u||\underline{\partial}_a\Gamma^{I_2}\partial^2u|\ \!\|_{L^2_f(\mathcal{H}_s)}\nonumber\\
\!\!\!&\lesssim&\!\!\!\sum_{\substack{|I_1|+|I_2|\le|I|+2\\a\in\{1,2\}}}\|(s/t)|\partial\Gamma^{I_1} u||\partial\Gamma^{I_2}u|+|\underline{\partial}_a\Gamma^{I_1} u||\partial\Gamma^{I_2}u|\ \!\|_{L^2_f(\mathcal{H}_s)}\nonumber\\
\!\!\!&\lesssim&\!\!\!\sum_{\substack{|I_1|\le N-7,|I_2|\le N\\a\in\{1,2\}}}\|s|\partial\Gamma^{I_1} u|+t|\underline{\partial}_a\Gamma^{I_1} u|\ \!\|_{L^\infty(\mathcal{H}_s)}\cdot\|(s/t)|\partial\Gamma^{I_2}u|\ \!\|_{L^2_f(\mathcal{H}_s)}\cdot s^{-1}\nonumber\\
\!\!\!&+&\!\!\!\sum_{\substack{|I_2|\le N-7,|I_1|\le N\\a\in\{1,2\}}}\|\ \!|\underline{\partial}_a\Gamma^{I_1} u|\ \!\|_{L^2_f(\mathcal{H}_s)}\cdot\|s|\partial\Gamma^{I_2}u|\ \!\|_{L^\infty(\mathcal{H}_s)}\cdot s^{-1}\lesssim (C_1\ez)^2s^{-1+\dz}\lesssim(C_1\ez)^2.
\een
Lemmas \ref{lemHo} and \ref{lemnf} yield
\beq\label{gzN_2uuL2}
\begin{split}
&\ \|\Gamma^I(N_2(u,u))\|_{L^2_f(\mathcal{H}_s)}\lesssim\sum_{I_1+I_2\le I}\|N_{2,I;I_1,I_2}(\Gamma^{I_1}u,\Gamma^{I_2}u)\|_{L^2_f(\mathcal{H}_s)}\\
\lesssim&\ \sum_{\substack{|I_1|+|I_2|\le|I|\\|J_1|=1,|J_2|\le 1}}\|(s^2/t^2)|\partial\Gamma^{I_1}u||\partial^2\Gamma^{I_2}u|+t^{-1}|\Gamma^{J_1}\Gamma^{I_1}u||\partial\Gamma^{J_2}\Gamma^{I_2}u|\ \!\|_{L^2_f(\mathcal{H}_s)}\\
\lesssim&\ \sum_{\substack{|I_1|+|I_2|\le|I|+1\\|J|=1}}\|(s^2/t^2)|\partial\Gamma^{I_1}u||\partial\Gamma^{I_2}u|+t^{-1}|\Gamma^{J}\Gamma^{I_1}u||\partial\Gamma^{I_2}u|\ \!\|_{L^2_f(\mathcal{H}_s)}\\
\lesssim&\ \sum_{\substack{|I_1|\le N-7,|I_2|\le N\\|J|=1}}\|s|\partial\Gamma^{I_1}u|+|\Gamma^J\Gamma^{I_1}u|\ \!\|_{L^\infty(\mathcal{H}_s)}\cdot\|(s/t)|\partial\Gamma^{I_2}u|\ \!\|_{L^2_f(\mathcal{H}_s)}\cdot s^{-1}\\
+&\ \sum_{\substack{|I_2|\le N-7,|I_1|\le N\\|J|=1}}\|t^{-1}|\Gamma^J\Gamma^{I_1}u|\ \!\|_{L^2_f(\mathcal{H}_s)}\cdot\|s|\partial\Gamma^{I_2}u|\ \!\|_{L^\infty(\mathcal{H}_s)}\cdot s^{-1}\lesssim\ (C_1\ez)^2s^{-1+\dz}\lesssim(C_1\ez)^2,
\end{split}
\eeq
where for any sufficiently smooth functions $w,z$, $N_{2,I;I_1,I_2}(w,z)=P_{\!2,I;I_1,I_2}^{\gz\az\bz}\partial_\gz w\partial_\az\partial_\bz z$ with $P_{\!2,I;I_1,I_2}^{\gz\az\bz}$ satisfying the null condition. Combining the above estimates, we obtain
\beq\label{E_1v-tv}
[\mathcal{E}_1(\Gamma^I(v-\tilde{v}),s)]^{1/2}\lesssim (C_1\ez)^2,\quad|I|\le N-5.
\eeq
This together with \eqref{E_0u-tu} and \eqref{imtutv} yields
\beq\label{imLuv}
[\mathcal{E}_0(\Gamma^Iu,s)]^{1/2}+[\mathcal{E}_1(\Gamma^Iv,s)]^{1/2}\lesssim\ez+(C_1\ez)^2,\quad |I|\le N-5.
\eeq
Hence we have strictly improved the estimates of $\sum_{|I|\le N-5}\l\{[\mathcal{E}_0(\Gamma^Iu,s)]^{1/2}+[\mathcal{E}_1(\Gamma^Iv,s)]^{1/2}\r\}$ in \eqref{Bsuv} if we choose $C_1>1$ sufficiently large and $0<\ez\ll C_1^{-1}$ sufficiently small. In addition, by Proposition \ref{sC: DL5.1con}, \eqref{'btv+tuR} and \eqref{EsQ1}-\eqref{EsQ3}, for $|I|\le N-5$, we obtain
\be
[\mathcal{E}_{con}(\Gamma^I\tilde{u},s)]^{1/2}\lesssim[\mathcal{E}_{con}(\Gamma^I\tilde{u},s_0)]^{1/2}+\int_{s_0}^s\tau\|(-\Box)\Gamma^I\tilde{u}\|_{L^2_f(\mathcal{H}_\tau)}{\rm{d}}\tau
\lesssim\ez+(C_1\ez)^2s^\dz.
\ee
This together with Proposition \ref{sC: DL5.2} implies that for $|I|\le N-5$
\beq\label{EtuL2}
\|(s/t)\Gamma^I\tilde{u}\|_{L^2_f(\mathcal{H}_s)}\lesssim\|(s_0/t)\Gamma^I\tilde{u}\|_{L^2_f(\mathcal{H}_{s_0})}+\int_{s_0}^s\f{[\mathcal{E}_{con}(\Gamma^I\tilde{u},\tau)]^{1/2}}{\tau}{\rm{d}}\tau
\lesssim\big(\ez+(C_1\ez)^2\big)s^\dz.
\eeq
Recall that $\tilde{u}-u=\f{1}{4}N_1(v,v)+N_2(u,v)$. We compute
\beq\label{gzu-tuL2}
\begin{split}
&\ \|\ \!|\Gamma^I(N_1(v,v))|+|\Gamma^I(N_2(u,v))|\ \!\|_{L^2_f(\mathcal{H}_s)}\lesssim\sum_{|I_1|+|I_2|\le|I|+2}\|\big(|\partial\Gamma^{I_1}v|+|\partial\Gamma^{I_1}u|\big)|\Gamma^{I_2}v|\ \!\|_{L^2_f(\mathcal{H}_s)}\\
\lesssim&\ \sum_{|I_1|\le N-7,|I_2|\le N}\|s\big(|\partial\Gamma^{I_1}v|+|\partial\Gamma^{I_1}u|\big)\|_{L^\infty(\mathcal{H}_s)}\cdot\|\ \!|\Gamma^{I_2}v|\ \!\|_{L^2_f(\mathcal{H}_s)}\cdot s^{-1}\\
+&\sum_{|I_2|\le N-7,|I_1|\le N}\|(s/t)\big(|\partial\Gamma^{I_1}v|+|\partial\Gamma^{I_1}u|\big)\|_{L^2_f(\mathcal{H}_s)}\cdot\|t|\Gamma^{I_2}v|\ \!\|_{L^\infty(\mathcal{H}_s)}\cdot s^{-1}\\
\lesssim&\ (C_1\ez)^2s^{-1+\dz}\lesssim(C_1\ez)^2,\quad |I|\le N-5.
\end{split}
\eeq
It follows from \eqref{gzu-tuL2} and \eqref{EtuL2} that
\beq\label{EuL2}
\|(s/t)\Gamma^Iu\|_{L^2_f(\mathcal{H}_s)}\lesssim\big(\ez+(C_1\ez)^2\big)s^\dz,\quad |I|\le N-5.
\eeq
By Lemma \ref{sC: DL2.4} and \eqref{EuL2}, we arrive at
\beq\label{EuLinf}
\|\Gamma^Iu\|_{L^\infty(\mathcal{H}_s)}\lesssim\big(\ez+(C_1\ez)^2\big)s^{-1+\dz},\quad |I|\le N-7.
\eeq

Combining \eqref{imTuv} and \eqref{imLuv}, we have strictly improved the bootstrap estimate \eqref{Bsuv} (here we choose $C_1>1$ sufficiently large and $0<\ez\ll C_1^{-1}$ sufficiently small). Hence the proof of Proposition \ref{propim} is completed. In addition, the estimates \eqref{L2Linf} and \eqref{EuLinf} hold for all $s\in[2,+\infty)$.

\section{Proof of the scattering result}\label{sSca}

In this section we prove the scattering result in Theorem \ref{thm1}.

Let $(u,v)$ be the solution to the following Cauchy problem in $\mathbb{R}^{1+2}$
\be
\l\{\begin{array}{rcl}
-\Box u&=&f_u,\\
-\Box v+v&=&f_v,
\end{array}\r.\quad\quad\quad\quad(u,\partial_tu,v,\partial_tv)|_{t=t_0=2}=(u^0,u^1,v^0,v^1)
\ee
with the initial data $(u^0,u^1,v^0,v^1)$ supported in the ball $\{x: |x|<1\}$. We denote 
\beq\label{s5: vecu}
\vec{u}=(u,\partial_tu)'=\l(\begin{array}{c}
u\\
\partial_tu
\end{array}\r),\quad\quad\vec{v}=(v,\partial_tv)'=\l(\begin{array}{c}
v\\
\partial_tv
\end{array}\r),  
\eeq
where $(a_1,a_2)'$ denotes the transpose of a vector $\vec{a}=(a_1,a_2)$ in $\mathbb{R}^2$. We also set
\beq\label{s5: vecF}
\vec{f}_u=(0,f_u)',\quad\quad\vec{f}_v=(0,f_v)',\quad\quad\vec{u}^0=(u^0,u^1)',\quad\quad\vec{v}^0=(v^0,v^1)'.
\eeq 
By the linear theory of wave and Klein-Gordon equations, we can write
\beq\label{s5: vecu=}
\vec{u}=\mathcal{S}(t-2)\vec{u}^0+\int_2^t\mathcal{S}(t-\tau)\vec{f}_u(\tau){\rm{d}}\tau,
\eeq
\beq\label{s5: vecv=}
\vec{v}=\mathcal{\tilde{S}}(t-2)\vec{v}^0+\int_2^t\mathcal{\tilde{S}}(t-\tau)\vec{f}_v(\tau){\rm{d}}\tau,
\eeq
where
\be
\mathcal{S}(t)=\l(\begin{array}{cc}
\cos(t\sqrt{-\Delta})&\f{\sin(t\sqrt{-\Delta})}{\sqrt{-\Delta}}\\
-\sqrt{-\Delta}\sin(t\sqrt{-\Delta})&\cos(t\sqrt{-\Delta})
\end{array}\r),\quad\quad\mathcal{\tilde{S}}(t)=\l(\begin{array}{cc}
\cos(t\langle\nabla\rangle)&\f{\sin(t\langle\nabla\rangle)}{\langle\nabla\rangle}\\
-\langle\nabla\rangle\sin(t\langle\nabla\rangle)&\cos(t\langle\nabla\rangle)
\end{array}\r).
\ee

Let $l\in\mathbb{N}$. We denote ${\bf H}^l(\mathbb{R}^2):=H^{l+1}(\mathbb{R}^2)\times H^l(\mathbb{R}^2)$ and ${\bf{H}}_{l}(\mathbb{R}^2):=\big(\dot{H}^{l+1}(\mathbb{R}^2)\cap\dot{H}^1(\mathbb{R}^2)\big)\times H^l(\mathbb{R}^2)$, where $H^k(\mathbb{R}^2), \dot{H}^{k}(\mathbb{R}^2), k\in\mathbb{N}$ denote the Sobolev spaces and homogeneous Sobolev spaces respectively.

\begin{lem}\label{lemSca}
The following statements hold:
\begin{itemize}
\item[i)] Let $l\in\mathbb{N}$ and $\dz>0$. For any $\mathbb{R}^2$-valued function $\vec{f}(\tau,x)=(f_1,f_2)'$ which is defined in $[2,+\infty)\times\mathbb{R}^2$ with support in $\mathcal{K}$ and satisfies $\vec{f}(\tau,\cdot)\in{\bf{H}}^l(\mathbb{R}^2)$ for any fixed $\tau\in[2,+\infty)$, any $t\in[2,+\infty)$, and any $4\le T_1<T_2<+\infty$, we have
\ben
&&\l\|\int_{T_1}^{T_2}\mathcal{S}(t-\tau)\vec{f}(\tau){\rm{d}}\tau\r\|_{{\bf{H}}_{l}(\mathbb{R}^2)}\lesssim T_1^{-\f{\dz}{2}}\sum_{k=0}^l\l(\int_{T_1^{\f{1}{2}}}^{T_2}\|\ \!|\vec{\nabla}^k\vec{f}|\ \!\|_{L^2_f(\mathcal{H}_s)}^2\cdot s^{1+2\dz}{\rm{d}}s\r)^{1/2},\\
&&\l\|\int_{T_1}^{T_2}\mathcal{\tilde{S}}(t-\tau)\vec{f}(\tau){\rm{d}}\tau\r\|_{{\bf{H}}^l(\mathbb{R}^2)}
\lesssim T_1^{-\f{\dz}{2}}\sum_{k=0}^l\l(\int_{T_1^{\f{1}{2}}}^{T_2}\|\ \!|\vec{\nabla}^k\vec{f}|+|f_1|\ \!\|_{L^2_f(\mathcal{H}_s)}^2\cdot s^{1+2\dz}{\rm{d}}s\r)^{1/2},
\een
where $|\vec{\nabla}^k\vec{f}|:=|\nabla^{k+1}f_1|+|\nabla^kf_2|$ and $\nabla=(\partial_1,\partial_2)$.

\item[ii)] Let $l\in\mathbb{N}$ and $\vec{u},\vec{v},\vec{f}_u,\vec{f}_v,\vec{u}^0,\vec{v}^0$ be as in \eqref{s5: vecu}-\eqref{s5: vecv=} with $\vec{u}^0\in{\bf{H}}_{l}(\mathbb{R}^2), \vec{v}^0\in{\bf{H}}^{l}(\mathbb{R}^2)$ supported in $\{x: |x|<1\}$ and $f_u,f_v$ supported in $\mathcal{K}$ satisfying $f_u(\tau),f_v(\tau)\in H^l(\mathbb{R}^2)$ for any fixed $\tau\in[2,+\infty)$. If for some $\dz>0$, it holds that
\be
M:=\int_2^4\|f(\tau,\cdot)\|_{L^2_x(\mathbb{R}^2)}{\rm{d}}\tau+\l(\int_{2}^{+\infty}\|f\|_{L^2_f(\mathcal{H}_s)}^2\cdot s^{1+2\dz}{\rm{d}}s\r)^{1/2}<+\infty,
\ee
where $f:=\sum_{k=0}^{l}\big(|\nabla^{k}f_u|+|\nabla^kf_v|\big)$, then the solution $(\vec{u},\vec{v})$ scatters to a free solution in ${\bf{H}}_{l}(\mathbb{R}^2)\times {\bf{H}}^{l}(\mathbb{R}^2)$, i.e., there exist $\vec{u}^*_0=(u^*_0,u^*_1)'\in{\bf{H}}_{l}(\mathbb{R}^2)$ and $\vec{v}^*_0=(v^*_0,v^*_1)'\in{\bf{H}}^{l}(\mathbb{R}^2)$ such that
\ben
\lim_{t\to+\infty}\|\vec{u}-\vec{u}^*\|_{{\bf{H}}_{l}(\mathbb{R}^2)}=0,\quad\quad\lim_{t\to+\infty}\|\vec{v}-\vec{v}^*\|_{{\bf{H}}^{l}(\mathbb{R}^2)}=0,
\een
where $\vec{u}^*=(u^*,\partial_tu^*)'$, $\vec{v}^*=(v^*,\partial_tv^*)'$, and $(u^*,v^*)$ is the solution to the $2D$ linear homogeneous wave-Klein-Gordon system with the initial data $(u^*_0,u^*_1,v^*_0,v^*_1)$ (prescribed on the time slice $t=t_0=2$).
\end{itemize}

\begin{proof}
$i)$ We only need to consider the case $l=0$. For any fixed $t\in[2,+\infty)$, let $\vec{U}(\tau):=\mathcal{S}(t-\tau)\vec{f}(\tau)=(U_1,U_2)'$. For any $4\le T_1<T_2<+\infty$, by standard energy inequalities, we have
\ben
\l\|\int_{T_1}^{T_2}\vec{U}(\tau){\rm{d}}\tau\r\|_{{\bf{H}}_0(\mathbb{R}^2)}&\lesssim&\l(\int_{\mathbb{R}^2}\Bigg\{\l|\int_{T_1}^{T_2}\nabla U_1(\tau){\rm{d}}\tau\r|^2+\l|\int_{T_1}^{T_2}U_2(\tau){\rm{d}}\tau\r|^2\Bigg\}{\rm{d}}x\r)^{1/2}\\
&\lesssim&\l(\int_{\mathbb{R}^2}\l\{\int_{T_1}^{T_2}|\nabla U_1(\tau)\cdot \tau^{\f{1+\dz}{2}}|^2{\rm{d}}\tau+\int_{T_1}^{T_2}|U_2(\tau)\cdot\tau^{\f{1+\dz}{2}}|^2{\rm{d}}\tau\r\}\cdot\l(\int_{T_1}^{T_2}\tau^{-(1+\dz)}{\rm{d}}\tau\r){\rm{d}}x\r)^{1/2}\\
&\lesssim&T_1^{-\f{\dz}{2}}\l(\int_{T_1}^{T_2}\int_{\mathbb{R}^2}\l\{|\nabla U_1(\tau)|^2+|U_2(\tau)|^2\r\}\cdot \tau^{1+\dz}{\rm{d}}x{\rm{d}}\tau\r)^{1/2}\\
&\lesssim&II:=T_1^{-\f{\dz}{2}}\l(\int_{T_1}^{T_2}\int_{r<\tau-1}\l\{|\nabla f_1(\tau)|^2+|f_2(\tau)|^2\r\}\cdot \tau^{1+\dz}{\rm{d}}x{\rm{d}}\tau\r)^{1/2}.
\een
By a change of variables $(\tau,x)\to(s,x)$ with $s=\sqrt{\tau^2-|x|^2}$, we obtain
\beq\label{EII}
\begin{split}
II\lesssim&\ T_1^{-\f{\dz}{2}}\l(\int_{T_1^{\f{1}{2}}}^{T_2}\int_{r<\f{s^2-1}{2}}\l\{|\nabla f_1|^2+|f_2|^2\r\}(\sqrt{s^2+|x|^2},x)\cdot \tau^{1+\dz}\f{s}{\tau}{\rm{d}}x{\rm{d}}s\r)^{1/2}\\
\lesssim&\ T_1^{-\f{\dz}{2}}\l(\int_{T_1^{\f{1}{2}}}^{T_2}\|\ \!|\nabla f_1|+|f_2|\ \!\|_{L^2_f(\mathcal{H}_s)}^2\cdot s^{1+2\dz}{\rm{d}}s\r)^{1/2}.
\end{split}
\eeq
Similarly,
\be
\l\|\int_{T_1}^{T_2}\mathcal{\tilde{S}}(t-\tau)\vec{f}(\tau){\rm{d}}\tau\r\|_{{\bf{H}}^0(\mathbb{R}^2)}
\lesssim T_1^{-\f{\dz}{2}}\l(\int_{T_1^{\f{1}{2}}}^{T_2}\|\ \!|\nabla f_1|+|f_2|+|f_1|\ \!\|_{L^2_f(\mathcal{H}_s)}^2\cdot s^{1+2\dz}{\rm{d}}s\r)^{1/2}.
\ee
$ii)$ Let $\vec{u}$, $\vec{v}$, $\vec{u}^0$, $\vec{v}^0$, $\vec{f}_u$, $\vec{f}_v$ be as in \eqref{s5: vecu}-\eqref{s5: vecv=} and 
\ben
\vec{u}^*_0:&=&\vec{u}^0+\int_2^{+\infty}\mathcal{S}(2-\tau)\vec{f}_u(\tau){\rm{d}}\tau,\quad\quad \vec{u}^*=\mathcal{S}(t-2)\vec{u}^*_0,\\
\vec{v}^*_0:&=&\vec{v}^0+\int_2^{+\infty}\mathcal{\tilde{S}}(2-\tau)\vec{f}_v(\tau){\rm{d}}\tau,\quad\quad \vec{v}^*=\mathcal{\tilde{S}}(t-2)\vec{v}^*_0.
\een
For any $4\le T_1<T_2<+\infty$, by $i)$, we have
\ben
&&\l\|\int_{T_1}^{T_2}\mathcal{S}(2-\tau)\vec{f}_u(\tau){\rm{d}}\tau\r\|_{{\bf{H}}_{l}(\mathbb{R}^2)}+\l\|\int_{T_1}^{T_2}\mathcal{\tilde{S}}(2-\tau)\vec{f}_v(\tau){\rm{d}}\tau\r\|_{{\bf{H}}^{l}(\mathbb{R}^2)}\\
&\lesssim& T_1^{-\f{\dz}{2}}\sum_{k=0}^{l}\l(\int_{T_1^{\f{1}{2}}}^{T_2}\|\ \!|\nabla^kf_u|+|\nabla^kf_v|\ \!\|_{L^2_f(\mathcal{H}_s)}^2\cdot s^{1+2\dz}{\rm{d}}s\r)^{1/2}\lesssim MT_1^{-\f{\dz}{2}}\to 0\quad\mathrm{as}\ \  T_1\to+\infty.
\een
Hence, $\vec{u}^*_0$ and $\vec{v}^*_0$ are well-defined in ${\bf{H}}_{l}(\mathbb{R}^2)$ and ${\bf{H}}^{l}(\mathbb{R}^2)$ respectively. Similarly, using $i)$ again, we have
\ben
&&\|\vec{u}-\vec{u}^*\|_{{\bf{H}}_{l}(\mathbb{R}^2)}=\l\|\int_t^{+\infty}\mathcal{S}(t-\tau)\vec{f}_u(\tau){\rm{d}}\tau\r\|_{{\bf{H}}_{l}(\mathbb{R}^2)}\lesssim Mt^{-\f{\dz}{2}}\to 0\quad\mathrm{as}\ \ t\to+\infty,\\
&&\|\vec{v}-\vec{v}^*\|_{{\bf{H}}^{l}(\mathbb{R}^2)}=\l\|\int_t^{+\infty}\mathcal{\tilde{S}}(t-\tau)\vec{f}_v(\tau){\rm{d}}\tau\r\|_{{\bf{H}}^{l}(\mathbb{R}^2)}\lesssim Mt^{-\f{\dz}{2}}\to 0\quad\mathrm{as}\ \ t\to+\infty.
\een
The proof is done.
\end{proof}
\end{lem}

$\it{Proof\;of\;the\;scattering\;result\;in\;Theorem\;1.}$

For $l\in\mathbb{N}$ we denote
\be
\mathcal{X}_{l}(\mathbb{R}^2):={\bf{H}}_{l}(\mathbb{R}^2)\times {\bf{H}}^{l}(\mathbb{R}^2)=\big(\dot{H}^{l+1}(\mathbb{R}^2)\cap\dot{H}^1(\mathbb{R}^2)\big)\times H^l(\mathbb{R}^2)\times H^{l+1}(\mathbb{R}^2)\times H^l(\mathbb{R}^2).
\ee
Let $(u,v)$ be the global solution to \eqref{qEquv}-\eqref{qini} given by Theorem \ref{thm1}. Denote $F_{\tilde{u}}:=(-\Box)\tilde{u}$ and $F_{\tilde{v}}:=(-\Box+1)\tilde{v}$, where $\tilde{u},\tilde{v}$ are as in Lemma \ref{lemNLuv}. Then by \eqref{'btv+tuR} and \eqref{EsQ1}-\eqref{EsQ3}, for any $s\in[2,+\infty)$ we have
\beq\label{EsgztF}
\sum_{|I|\le N-5}\|\ \!|\Gamma^IF_{\tilde{u}}|+|\Gamma^IF_{\tilde{v}}|\ \!\|_{L^2_f(\mathcal{H}_s)}\lesssim s^{-2+\dz}.
\eeq
It follows from \eqref{EsgztF} and Lemma \ref{lemSca} that $(\tilde{u},\partial_t\tilde{u},\tilde{v},\partial_t\tilde{v})$ scatters to a free solution in $\mathcal{X}_{N-5}(\mathbb{R}^2)$. It remains to show that
\beq\label{u-tuXn}
\lim_{t\to+\infty}\|(u,\partial_tu,v,\partial_tv)-(\tilde{u},\partial_t\tilde{u},\tilde{v},\partial_t\tilde{v})\|_{\mathcal{X}_{N-5}(\mathbb{R}^2)}=0.
\eeq
By the proof of \eqref{gzu-tuL2} and \eqref{gzN_2uuL2}, for any $s\in[2,+\infty)$, it holds that
\beq\label{gzFu+N_2u}
\|\ \!|\Gamma^IF_u|+|\Gamma^I(N_2(u,u))|\ \!\|_{L^2_f(\mathcal{H}_s)}\lesssim s^{-1+\dz},\quad |I|\le N-2.
\eeq
By \eqref{L2Linf}, for $|I|\le N-2$ and any $s\in[2,+\infty)$ we have
\beq\label{gzN_1vu}
\begin{split}
\!\!\!\!\!\!\|\Gamma^I(N_1(v,u))\|_{L^2_f(\mathcal{H}_s)}\lesssim&\ \sum_{|I_1|+|I_2|\le N}\|\ \!|\Gamma^{I_1}v||\partial\Gamma^{I_2}u|\ \!\|_{L^2_f(\mathcal{H}_s)}\\
\!\!\!\!\!\!\lesssim&\ \sum_{|I_1|\le N-7,|I_2|\le N}\|t|\Gamma^{I_1}v|\ \!\|_{L^\infty(\mathcal{H}_s)}\cdot\|(s/t)|\partial\Gamma^{I_2}u|\ \!\|_{L^2_f(\mathcal{H}_s)}\cdot s^{-1}\\
\!\!\!\!\!\!+&\ \sum_{|I_2|\le N-7,|I_1|\le N}\|\ \!|\Gamma^{I_1}v|\ \!\|_{L^2_f(\mathcal{H}_s)}\cdot\|s|\partial\Gamma^{I_2}u|\ \!\|_{L^\infty(\mathcal{H}_s)}\cdot s^{-1}\lesssim s^{-1+\dz}.
\end{split}
\eeq
Combining \eqref{gzFu+N_2u} and \eqref{gzN_1vu}, we obtain that for all $s\in[2,+\infty)$
\beq\label{gzLF_uv}
\sum_{|I|\le N-2}\|\ \!|\Gamma^IF_u|+|\Gamma^IF_v|\ \!\|_{L^2_f(\mathcal{H}_s)}\lesssim s^{-1+\dz}.
\eeq
Denote $G:=\sum_{|I|\le N-2}(|\Gamma^IF_u|+|\Gamma^IF_v|)$. Then for any $t\in[4,+\infty)$, by a change of variables $(\tau,x)\to(s,x)$ with $s=\sqrt{\tau^2-|x|^2}$ (similar to the proof of \eqref{EII}), we have
\ben
&&\int_4^t\|G(\tau)\|_{L^2_x(\mathbb{R}^2)}\cdot\tau^{\f{1+\dz}{2}}\cdot\tau^{-\f{1+\dz}{2}}{\rm{d}}\tau\lesssim\l(\int_4^t\int_{r<\tau-1}|G(\tau,x)|^2\cdot\tau^{1+\dz}{\rm{d}}\tau\r)^{1/2}\nonumber\\
&\lesssim&\l(\int_2^t\|G\|_{L^2_f(\mathcal{H}_s)}^2\cdot s^{1+2\dz}{\rm{d}}s\r)^{1/2}\lesssim t^{2\dz},
\een
where we use \eqref{gzLF_uv}. Hence by the standard energy estimates for wave and Klein-Gordon equations, for any $t\in[4,+\infty)$ we arrive at
\beq\label{Esdktuv}
\sum_{|I|\le N-2}\|\big(|\partial\Gamma^Iu|+|\partial\Gamma^Iv|+|\Gamma^Iv|\big)(t)\|_{L^2_x(\mathbb{R}^2)}\lesssim\ez+\int_2^t\|G(\tau)\|_{L^2_x(\mathbb{R}^2)}{\rm{d}}\tau\lesssim  t^{2\dz}.
\eeq
By \eqref{L2Linf} and \eqref{Esdktuv}, we derive
\beq\label{u-tuXn1}
\begin{split}
&\ \sum_{|I|\le N-4}\|\big(|\Gamma^I(N_1(v,v))|+|\Gamma^I(N_2(u,v))|+|\Gamma^I(N_2(u,u))|\big)(t)\|_{L^2_x(\mathbb{R}^2)}\\
\lesssim&\ \sum_{|I_1|+|I_2|\le N-3}\|\big(|\partial\Gamma^{I_1}v|+|\partial\Gamma^{I_1}u|\big)\big(|\partial\Gamma^{I_2}v|+|\partial\Gamma^{I_2}u|\big)(t)\|_{L^2_x(\mathbb{R}^2)}\\
\lesssim&\ \sum_{|I_1|\le N-7,|I_2|\le N-3}\|\big(|\partial\Gamma^{I_1}v|+|\partial\Gamma^{I_1}u|\big)(t)\|_{L^\infty_x(\mathbb{R}^2)}\cdot\|\big(|\partial\Gamma^{I_2}v|+|\partial\Gamma^{I_2}u|\big)(t)\|_{L^2_x(\mathbb{R}^2)}\\
\lesssim&\ \ t^{-1/2+2\dz}.
\end{split}
\eeq
Recall that $\tilde{v}=v-N_2(u,u)$ and $\tilde{u}=u+\f{1}{4}N_1(v,v)+N_2(u,v)$ (see Lemma \ref{lemNLuv}). Hence \eqref{u-tuXn} follows from \eqref{u-tuXn1}, and we conclude that $(u,\partial_tu,v,\partial_tv)$ scatters to a free solution in $\mathcal{X}_{N-5}(\mathbb{R}^2)$.

Qian Zhang

School of Mathematical Sciences, Tiangong University, Tianjin, 300387, PR China


Email: m201011131071\_zq@126.com


\begin{thebibliography}{99}


\bibitem{Al1} \label{Al1} Alinhac, S.:
\textit{The null condition for quasilinear wave equations in two space dimensions I.}
Invent. Math. $\mathbf{145}$ (2001), no. 3, 597-618.

\bibitem{Al2} \label{Al2} Alinhac, S.:
\textit{The null condition for quasi linear wave equations in two space dimensions II.}
Amer. J. Math. $\mathbf{123}$ (2001), no. 6, 1071-1101.

\bibitem{Al02} \label{Al02} Alinhac, S.:
\textit{A minicourse on the global existence and blowup of classical solutions to multidimensional quasilinear wave equations.}
Journ\'ees "\'Equations aux D\'eriv\'ees Partielles" (Forges-les-Eaux, 2002), Exp. No. I, 33 pp., Univ. of Nantes, 2002.

\bibitem{Al03} \label{Al03} Alinhac, S.:
\textit{An example of blowup at infinity for a quasilinear wave equation.}
Autour de l'analyse microlocale, Ast\'erisque (2003), no. 284, 1-91.

\bibitem{Al10} \label{Al10} Alinhac, S.:
\textit{Geometric Analysis of Hyperbolic Differential Equations: An introduction.}
London Math. Soc. Lecture Note Ser., vol. 374, Cambridge University Press, Cambridge, 2010.



\bibitem{CLM} \label{CLM} Cai, Y., Lei, Z., Masmoudi, N.:
\textit{Global well-posedness for $2D$ nonlinear wave equations without compact support.}
J. Math. Pures Appl. (9) $\mathbf{114}$ (2018), 211-234.




\bibitem{Ch} \label{Ch} Christodoulou, D.:  
\textit{Global solutions of nonlinear hyperbolic equations for small initial data.} 
Comm. Pure Appl. Math. $\mathbf{39}$ (1986), no. 2, 267-282.


\bibitem{DFX} \label{DFX} Delort, J.-M., Fang, D., Xue, R.:
\textit{Global existence of small solutions for quadratic quasilinear Klein-Gordon systems in two space dimensions}
J. Funct. Anal. $\mathbf{211}$ (2004), no. 2, 288-323.


\bibitem{D21} \label{D21} Dong, S.:,
\textit{Global solution to the wave and Klein-Gordon system under null condition in dimension two.} 
J. Funct. Anal. $\mathbf{281}$ (2021), no.11, Paper No. 109232, 29 pp.















\bibitem{Do21} \label{Do21} Dong, S.: 
\textit{Asymptotic behavior of the solution to the Klein-Gordon-Zakharov model in dimension two.} 
Comm. Math. Phys. $\mathbf{384}$ (2021), no. 1, 587-607.


\bibitem{DLL} \label{DLL} Dong, S., LeFloch, P.G., Lei, Z.:
\textit{The top-order energy of quasilinear wave equations in two space dimensions is uniformly bounded.}
Preprint, arXiv:2103.07867.

\bibitem{DW20} \label{DW20} Dong, S., Wyatt, Z.: 
\textit{Stability of a coupled wave-Klein-Gordon system with quadratic nonlinearities.} 
J. Differential Equations $\mathbf{269}$ (2020), no. 9, 7470-7497.





\bibitem{DWd} \label{DWd} Dong, S., Wyatt, Z.: 
\textit{Hidden structure and sharp asymptotics for the Dirac-Klein-Gordon system in two space dimensions.} Preprint, arXiv:2105.13780. To appear in Ann. Inst. H. Poincar\'e C Anal. Non Lin\'eaire.

\bibitem{DuM} \label{DuM} Duan, S., Ma, Y.: 
\textit{Global solutions of wave-Klein-Gordon systems in $2+1$ dimensional space-time with strong couplings in divergence form.} 
SIAM J. Math. Anal. $\mathbf{54}$ (2022), no. 3, 2691-2726.

\bibitem{FWY} \label{FWY} Fang, A., Wang, Q., Yang, S.:
\textit{Global solution for massive Maxwell-Klein-Gordon equations with large Maxwell field.}
Ann. PDE $\mathbf{7}$ (2021), no. 1, Paper No. 3, 69 pp.

\bibitem{Geo} \label{Geo} Georgiev, V.:
\textit{Global solution of the system of wave and Klein-Gordon equations.}
Math. Z. $\mathbf{203}$ (1990), no. 4, 683-698.




\bibitem{Ho} \label{Ho} H\"ormander, L.:
\textit{Lectures on nonlinear hyperbolic differential equations.}
Math\'ematiques \& Applications (Berlin) [Mathematics \& Applications], vol. 26, Springer-Verlag, Berlin, 1997, viii+289 pp.


\bibitem{HY} \label{HY} Hou, F., Yin, H.:
\textit{Global small data smooth solutions of $2$-$D$ null-form wave equations with non-compactly supported initial data.}
J. Differential Equations $\mathbf{268}$ (2020), no. 2, 490-512.




\bibitem{IS} \label{IS} Ifrim, M., Stingo, A.:
\textit{Almost global well-posedness for quasilinear strongly coupled wave-Klein-Gordon systems in two space dimensions.}
Preprint, arXiv:1910.12673.



\bibitem{Ion} \label{Ion} Ionescu, A.D., Pausader, B.:
\textit{On the global regularity for a wave-Klein-Gordon coupled system.}
Acta Math. Sin. (Engl. Ser.) $\mathbf{35}$ (2019), no. 6, 933-986.




\bibitem{J81}\label{J81} John F.:
\textit{Blow up of solutions for quasi-linear wave equations in three space dimensions.} 
Comm. Pure Appl. Math. $\mathbf{34}$ (1981), no. 1, 29-51.


\bibitem{Kat} \label{Kat} Katayama, S.:
\textit{Global existence for coupled systems of nonlinear wave and Klein-Gordon equations in three space dimensions.}
Math. Z. $\mathbf{270}$ (2012), no. 1-2, 487-513.



\bibitem{K86} \label{K86} Klainerman, S.:
\textit{The null condition and global existence to nonlinear wave equations.}
Nonlinear systems of partial differential equations in applied mathematics, Part 1 (Santa Fe, N.M., 1984), 293-326, Lectures in Appl. Math., vol. 23, Amer. Math. Soc., Providence, RI, 1986.






\bibitem{K2} \label{K2} Klainerman, S.: 
\textit{Global existence of small amplitude solutions to nonlinear Klein-Gordon equations in four space-time dimensions.} 
Comm. Pure Appl. Math. $\mathbf{38}$ (1985), no. 5, 631-641.


\bibitem{K93} \label{K93} Klainerman, S.:
\textit{Remark on the asymptotic behavior of the Klein-Gordon equation in $\mathbb{R}^{n+1}$.}
Comm. Pure Appl. Math. $\mathbf{46}$ (1993), no. 2, 137-144.


\bibitem{KWY} \label{KWY} Klainerman, S., Wang, Q., Yang, S.: 
\textit{Global solution for massive Maxwell-Klein-Gordon equations.} 
Comm. Pure Appl. Math. $\mathbf{73}$ (2020), no. 1, 63-109.


\bibitem{LM14} \label{LM14} LeFloch, P.G., Ma, Y.:
\textit{The hyperboloidal foliation method.}
Series in Applied and Computational Mathematics, World Scientific Press, Hackensack, NJ (2014).


\bibitem{LM16} \label{LM16} LeFloch, P.G., Ma, Y.:
\textit{The global nonlinear stability of Minkowski space for self-gravitating massive fields. The Wave-Klein–Gordon Model.}
Comm. Math. Phys. $\mathbf{346}$ (2016), no. 2, 603-665.




\bibitem{Ld} \label{Ld} Li, D.:
\textit{Uniform estimates for $2D$ quasilinear wave.}
Adv. Math. $\mathbf{428}$ (2023), Paper No. 109157, 41 pp.


\bibitem{Lind} \label{Lind} Lindblad, H.: 
\textit{Global solutions of quasilinear wave equations.}
Am. J. Math. $\mathbf{130}$ (2008), no. 1, 115-157. 


\bibitem{M17} \label{M17} Ma, Y.:  
\textit{Global solutions of quasilinear wave-Klein-Gordon system in two-space dimension: technical tools,}
J. Hyperbolic Differ. Equ. $\mathbf{14}$ (4) (2017), 591-625.




\bibitem{M17'} \label{M17'} Ma, Y.:
\textit{Global solutions of quasilinear wave-Klein-Gordon system in two-space dimension: completion of the proof.}
J. Hyperbolic Differ. Equ. $\mathbf{14}$ (2017), no. 4, 627-670.



\bibitem{M19} \label{M19} Ma, Y.: 
\textit{Global solutions of nonlinear wave-Klein-Gordon system in two spatial dimensions: weak coupling case.}
Preprint, arXiv:1907.03516.

\bibitem{M21} \label{M21} Ma, Y.: 
\textit{Global solutions of nonlinear wave-Klein-Gordon system in two spatial dimensions: a prototype of strong coupling case.}
J. Differential Equations $\mathbf{287}$ (2021), 236-294.


\bibitem{OTT} \label{OTT} Ozawa, T., Tsutaya, K.; Tsutsumi, Y.:
\textit{Normal form and global solutions for the Klein-Gordon-Zakharov equations.} 
Ann. Inst. H. Poincar\'e C Anal. Non Lin\'eaire $\mathbf{12}$ (1995), no. 4, 459-503.


\bibitem{St18} \label{St18} Stingo, A.:
\textit{Global existence of small amplitude solutions for a model quadratic quasi-linear coupled wave-Klein-Gordon system in two space dimension, with mildly decaying Cauchy data.}
Preprint, arXiv:1810.10235. To appear in Mem. Amer. Math. Soc.


\bibitem{Wa} \label{Wa} Wang, Q.:
\textit{An intrinsic hyperboloid approach for Einstein Klein-Gordon equations.}
J. Differential Geom. $\mathbf{115}$ (2020), no. 1, 27-109.

\bibitem{Wo} \label{Wo} Wong, W.:
\textit{Small data global existence and decay for two dimensional wave maps.}
Preprint, arXiv:1712.07684. To appear in Ann. H. Lebesgue

\bibitem{Zha} \label{Zha} Zha, D.:
\textit{Global and almost global existence for general quasilinear wave equations in two space dimensions.}
 J. Math. Pures Appl. (9) $\mathbf{123}$ (2019), 270-299.

 
\end{thebibliography}
\end{document}